\documentclass[11pt]{article}

\usepackage{fullwidth}
\usepackage[top=1.5in, bottom=1.5in, left=1in, right=1in]{geometry}
\usepackage{xypic}
\usepackage{amsgen, amsmath, amstext, amsbsy, amsopn, amsfonts, amssymb, comment}
\usepackage{eepic}
\usepackage{epsf}
\xyoption{all}
\usepackage{tikz}
\usetikzlibrary{calc,through,backgrounds,shapes,matrix}

\usepackage[backref=page]{hyperref}
\usepackage{url, hypcap}
\hypersetup{colorlinks=true, citecolor=darkblue, linkcolor=darkblue}
\usepackage{mathrsfs}

\definecolor{darkblue}{rgb}{0.0,0,0.7} 
\definecolor{darkred}{rgb}{0.7,0,0} 
\newcommand{\darkred}{\color{darkred}} 
\definecolor{lightgrey}{rgb}{0.7,0.7,0.7} 

\newcommand{\defn}[1]{\emph{\darkred #1}} 

\def\Box{\square}

\def\mapright#1{\smash{\mathop{\longrightarrow}\limits^{#1}}}
\def\tra#1{\smash{\mathop{\mid\kern
-1pt\joinrel\relbar\joinrel\relbar}\limits^{*}_{#1}}}
\def\longtra#1{\smash{\mathop{\mid\kern
-1pt\joinrel\relbar\joinrel\relbar\joinrel\relbar}\limits^{*}_{#1}}}
\def\vlongtra#1{\smash{\mathop{\mid\kern
-1pt\joinrel\relbar\joinrel\relbar\joinrel\relbar\joinrel\relbar}\limits^{*}_{#1}}}
\def\vvlongtra#1{\smash{\mathop{\mid\kern
-1pt\joinrel\relbar\joinrel\relbar\joinrel\relbar\joinrel\relbar\joinrel\relbar}\limits^{*}_{#1}}}
\def\vvvlongtra#1{\smash{\mathop{\mid\kern
-1pt\joinrel\relbar\joinrel\relbar\joinrel\relbar\joinrel\relbar\joinrel\relbar\joinrel\relbar}\limits^{*}_{#1}}}
\def\etra#1{\smash{\mathop{\mid\kern
-1pt\joinrel\relbar\joinrel\relbar}\limits_{#1}}}

\def\iff{\Leftrightarrow}
\def\Rw{\Rightarrow}
\def\oo{\overline}

\def\I{{\cal{I}}}
\def\L{{\cal{L}}}

\def\N{\mathbb{N}}

\def\sem{\mbox{Sem}}
\def\res{\operatorname{Res}}
\def\rg{\operatorname{RG}}

\def\rc{\operatorname{RC}}

\def\src{\operatorname{SRC}}

\def\lcs{\mbox{lcs}}

\def\ker{\mathrm{ker}}

\def\cay{\mbox{Cay}}

\def\R{{\cal{R}}}

\def\D{{\bf D}}
\def\V{{\bf V}}

\def\p{\varphi}

\def\inv{^{-1}}

\def\bi{\begin{itemize}}
\def\ei{\end{itemize}}
\def\beq{\begin{equation}}
\def\eeq{\end{equation}}

\def\RZ{{\bf RZ}}
\def\RZnil{{\bf RZ}^{\bf N}}
\def\RLM{\mathrm{RLM}}
\def\Vert{\mathrm{Vert}}

\newtheorem{T}{Theorem}[section]
\newcommand{\bt}{\begin{T}}
\newcommand{\et}{\end{T}}
\newcommand{\ftd}{$\square$\end{T}}

\newtheorem{Proposition}[T]{Proposition}
\newcommand{\bp}{\begin{Proposition}}
\newcommand{\ep}{\end{Proposition}}
\newcommand{\fpd}{$\square$\end{Proposition}}

\newtheorem{Lemma}[T]{Lemma}
\newcommand{\bl}{\begin{Lemma}}
\newcommand{\el}{\end{Lemma}}
\newcommand{\fld}{$\square$\end{Lemma}}

\newtheorem{Corol}[T]{Corollary}
\newcommand{\bc}{\begin{Corol}}
\newcommand{\ec}{\end{Corol}}
\newcommand{\fcd}{$\square$\end{Corol}}

\newtheorem{Definition}[T]{Definition}

\newtheorem{Remark}[T]{Remark}
\newcommand{\br}{\begin{Remark}}
\newcommand{\er}{\end{Remark}}
\newcommand{\frd}{$\square$\end{Remark}}

\newtheorem{Example}[T]{Example}
\newcommand{\be}{\begin{Example}}
\newcommand{\ee}{\end{Example}}

\newtheorem{Problem}[T]{Problem}
\newcommand{\bq}{\begin{Problem}}
\newcommand{\eq}{\end{Problem}}

\newcommand{\proof}
   {\par\medbreak\noindent{\bf Proof}.\enspace}

\newcommand{\qed}{
$\Box$
\par\bigbreak}

\normalbaselines
\def\abstract#1{\par\bigskip
\begingroup\small
\baselineskip=12truept
\begin{center}ABSTRACT\end{center}
\par\medskip\par\noindent
\null\hfill\hbox{\vbox{\hsize=5truein\noindent#1}}
\hfill\null\par\endgroup\par}

\numberwithin{figure}{section}
\numberwithin{equation}{section}

\title{Random walks on semaphore codes and delay de Bruijn semigroups}
\author{{\bf John Rhodes}\\ 
{\em Department of Mathematics, University of California, Berkeley,}\\ 
{\em CA 94720, U.S.A.}\\
{\em email:} rhodes@math.berkeley.edu, BlvdBastille@aol.com\\
$ $\\
{\bf Anne Schilling}\\
{\em Department of Mathematics, University of California, Davis,}\\
{\em One Shields Ave., Davis, CA 95616-8633, U.S.A.}\\
{\em email:} anne@math.ucdavis.edu\\
$ $\\
{\bf Pedro V. Silva}\\
{\em Centro de
Matem\'{a}tica, Faculdade de Ci\^{e}ncias, Universidade do
Porto,}\\ {\em R. Campo Alegre 687, 4169-007 Porto, Portugal}\\
{\em email:} pvsilva@fc.up.pt}

\date{\today}

\begin{document}

\maketitle

\abstract{
We develop a new approach to random walks on de Bruijn graphs over the alphabet $A$ through right congruences 
on $A^k$, defined using the natural right action of $A^+$. A major role is played by special right congruences, which 
correspond to semaphore codes and allow an easier computation of the hitting time. We show how right 
congruences can be approximated by special right congruences.
}

\section{Introduction}

In graph theory, a $k$-dimensional \defn{de Bruijn graph} over the alphabet $A$ is a directed graph representing 
overlaps between sequences of symbols~\cite{deBruijn:1946, Good:1946}. The de Bruijn graph has $|A|^k$ vertices, 
given by all words of length $k$ in the alphabet $A$. There is an edge from vertex $a_1 \ldots a_k \in A^k$ to vertex 
$a_2 \ldots a_k a \in A^k$ for every $a \in A$. An important question for cryptography and networking is that
of de Bruijn sequences. A de Bruijn sequence is a cyclic word of length $|A|^k$  such that every possible word of length
$k$ over the alphabet $A$ appears once and exactly once (see~\cite{Ralston:1982} for a review on
de Bruijn sequences). Obviously, a de Bruijn sequence corresponds to a Eulerian path in the de Bruijn graph.

Here we are interested in random walks on the de Bruijn graph $\Gamma$. To an edge $v \overset{a}{\longrightarrow} w$
in $\Gamma$ we associate a probability $0 \le \pi(a) \le 1$, satisfying $\sum_{a \in A} \pi(a)=1$. This gives rise 
to the \defn{de Bruijn--Bernoulli process} (see for example~\cite{AS:2013,ABCN:2015}): if we are at vertex $v$ at a given time,
then with probability $\pi(a)$ we go to vertex $w$ where $v \overset{a}{\longrightarrow} w$ is an edge in $\Gamma$.
The \defn{transition matrix} $\mathcal{T} = (\mathcal{T}_{v,w})_{v,w \in A^k}$ encodes the transition probabilities,
that is, $\mathcal{T}_{v,w} = \pi(a)$ if $v \overset{a}{\longrightarrow} w$. Given a random walk, an important question
is to determine the \defn{stationary distribution}, which intuitively is the state that is reached after taking many
steps in the random walk. Mathematically, the stationary distribution is the vector $I$ such that
$I \mathcal{T} = I$. In other words, $I$ is the left eigenvector of $\mathcal{T}$ with eigenvalue one.
In the case of the de Bruijn--Bernoulli random walk, the stationary distribution $I \in A^k$ is multiplicative~\cite{AS:2013}
\[
	I = \Bigl( \prod_{a\in w} \pi(a) \Bigr)_{w \in A^k}.
\]

We can reformulate the random walk on the de Bruijn graph in algebraic terms. Namely, let us define the right
action of $A$ on $A^k$ by
\[
	a_1 \ldots a_k . a = a_2 \ldots a_k a
\]
for $a_1 \ldots a_k \in A^k$ and $a\in A$. This induces the action of the semigroup
$F(|A|,k) := A^1 \cup A^2 \cup \cdots \cup A^k = A^{\le k}$ of all words in $A$ of length $1,2,\ldots,k$ with
the multiplication $\cdot$ being concatenation and taking the last $k$ letters if the length is bigger than $k$.
For example, if $A=\{a,b\}$ and $k=3$, we have $ab \cdot ba = bba$ in $F(2,3)$. In this formulation, it is clear
that the walk in $j$ steps given by $a_1 \cdots a_j$ acts as a constant map (i.e., is independent of
the initial vertex) if and only if $j=k$. We call such elements \defn{resets}.

\bigskip

Random walks on de Bruijn graphs are a ``classical'' subject. However, in applications it is
\defn{right congruences}\footnote{An equivalence relation is a right congruence if it preserves the right action of a semigroup.
See Definition~\ref{definition.right congruence} for more details.}~\cite{A:1985,ML:1993,ML:1997,RS:2008} 
on $A^k$ (denoted by $\rc(A^k)$) under the faithful action of $F(|A|,k)$ and the associated random walks on their 
congruence classes that are important. Intuitively, these are the finite semigroups for which any product of $k$ 
elements act like constant maps on $A^k$, but because of the right congruence some products of length less than $k$
might be constant. Right congruences are a standard idea in finite state machines or finite automata 
theory~\cite{Rhodes:2010}. In finite state machines, they are used in passing to the unique minimal 
automata doing the same computation.
For example, assume one has a stream of data (e.g. chemical data on waste water being
emptied into a river).  Assume that there exist a positive integer $k$, so that only the $k$ most recent
symbols of data matter. Then there is a function $f \colon A^k \to D$, where $D$ is the data set.
The function could be of the form $f(a_1, \ldots, a_k)$ is ok or not ok (that is, $D$ is a two element set) 
depending on whether this recent $k$ long data meets EPA standards. Then the function $f$ gives an 
equivalence relation $\sim$ on $A^k$ given by $s\sim t$ if and only if $f(s) = f(t)$. In addition, there 
is a \textit{unique} maximal refinement of $\sim$ which is a right congruence (that is, the best lower 
approximation by a right congruence) $R$, namely $sRt$ for $s, t \in A^k$ if and only if
for all strings $u \in A^*$ we have $s.u \sim t.u$  or equivalently  $f(su) = f(tu)$.
Here . is the multiplication in $F(|A|, k)$. Then $(A^k/R, F(|A|,k))$ can compute the function $f$ since
$f$ factors through the $R$ classes (take $u$ to be 1). See~\cite{Rhodes:2010} for more details.

Consider the right congruence in $\rc(A^3)$ with $A=\{a,b\}$ defined by the congruence classes
\begin{equation}
\label{equation.rc example}
	\{aaa,baa,aba\}, \{bba\}, \{aab,bab\}, \{abb\}, \{bbb\}.
\end{equation}
\begin{figure}[t]
\begin{center}
\scalebox{0.7}{
\begin{tikzpicture}[>=latex,line join=bevel,]
\node (node_4) at (54.225bp,77.0bp) [draw,draw=none] {$\{bbb\}$};
  \node (node_3) at (31.225bp,151.0bp) [draw,draw=none] {$\{abb\}$};
  \node (node_2) at (31.225bp,223.0bp) [draw,draw=none] {$\{aab,bab\}$};
  \node (node_1) at (31.225bp,7.0bp) [draw,draw=none] {$\{bba\}$};
  \node (node_0) at (72.225bp,294.0bp) [draw,draw=none] {$\hspace{-1.5cm}\{aaa,baa,aba\}$};
  \draw [black,->] (node_2) ..controls (31.225bp,205.48bp) and (31.225bp,183.55bp)  .. (node_3);
  \definecolor{strokecol}{rgb}{0.0,0.0,0.0};
  \pgfsetstrokecolor{strokecol}
  \draw (40.225bp,188.0bp) node {$b$};
  \draw [black,->] (node_0) ..controls (56.683bp,288.75bp) and (39.561bp,281.17bp)  .. (32.225bp,268.0bp) .. controls (27.464bp,259.46bp) and (27.12bp,248.51bp)  .. (node_2);
  \draw (41.225bp,258.0bp) node {$b$};
  \draw [black,->] (node_3) ..controls (26.004bp,129.94bp) and (18.916bp,96.353bp)  .. (21.225bp,68.0bp) .. controls (22.452bp,52.923bp) and (25.359bp,35.956bp)  .. (node_1);
  \draw (30.225bp,77.0bp) node {$a$};
  \draw [black,->] (node_2) ..controls (40.939bp,234.62bp) and (46.248bp,241.49bp)  .. (50.225bp,248.0bp) .. controls (56.059bp,257.56bp) and (61.499bp,268.9bp)  .. (node_0);
  \draw (69.225bp,258.0bp) node {$a$};
  \draw [black,->] (node_1) ..controls (25.499bp,18.844bp) and (22.415bp,25.733bp)  .. (20.225bp,32.0bp) .. controls (9.3045bp,63.242bp) and (5.4821bp,71.18bp)  .. (1.2245bp,104.0bp) .. controls (-3.6735bp,141.76bp) and (12.884bp,184.79bp)  .. (node_2);
  \draw (10.225bp,114.0bp) node {$b$};
  \draw [black,->] (node_1) ..controls (58.017bp,13.525bp) and (124.22bp,33.35bp)  .. (124.22bp,77.0bp) .. controls (124.22bp,223.0bp) and (124.22bp,223.0bp)  .. (124.22bp,223.0bp) .. controls (124.22bp,248.54bp) and (102.05bp,271.01bp)  .. (node_0);
  \draw (133.22bp,151.0bp) node {$a$};
  \draw [black,->] (node_3) ..controls (36.803bp,133.05bp) and (44.171bp,109.35bp)  .. (node_4);
  \draw (54.225bp,114.0bp) node {$b$};
  \draw [black,->] (node_4) ..controls (48.507bp,59.599bp) and (41.54bp,38.395bp)  .. (node_1);
  \draw (53.225bp,41.0bp) node {$a$};
  \draw [black,->] (node_0) ..controls (86.459bp,302.57bp) and (96.225bp,301.03bp)  .. (96.225bp,294.0bp) .. controls (96.225bp,289.83bp) and (92.782bp,287.59bp)  .. (node_0);
  \draw (105.22bp,294.0bp) node {$a$};
  \draw [black,->] (node_4) ..controls (68.459bp,85.569bp) and (78.225bp,84.031bp)  .. (78.225bp,77.0bp) .. controls (78.225bp,72.825bp) and (74.782bp,70.587bp)  .. (node_4);
  \draw (87.225bp,77.0bp) node {$b$};
\end{tikzpicture}
}
\end{center}
\caption{The transition graph for the congruence of Equation~\eqref{equation.rc example}.
\label{figure.transition graph}}
\end{figure}
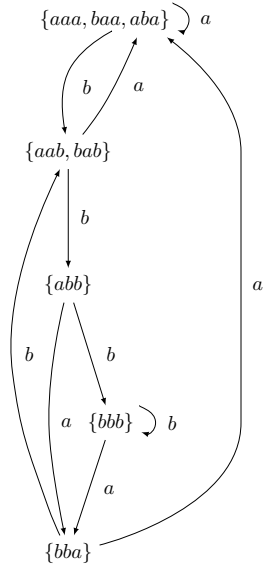
It is not hard to check that if $w, v \in A^3$ are in the same congruence class, then $w\cdot z$ and $v\cdot z$
for $z\in F(2,3)$ are also in the same congruence class, proving that~\eqref{equation.rc example} is indeed in $\rc(A^3)$.
The transition graph is given in Figure~\ref{figure.transition graph} and 
the transition matrix of the associated random walk is
\[
	\mathcal{T} = \begin{pmatrix} \pi(a)&0&\pi(b)&0&0\\
	\pi(a)&0&\pi(b)&0&0\\
	\pi(a)&0&0&\pi(b)&0\\
	0&\pi(a)&0&0&\pi(b)\\
	0&\pi(a)&0&0&\pi(b)
	 \end{pmatrix}\;.
\]
By lumping~\cite{KS:1976,LPW:2009}, we can obtain the stationary distribution for $\mathcal{T}$ from the 
stationary distribution of the de Bruijn--Bernoulli stationary distribution by adding the product distributions for
each member of a congruence class. In our example
\begin{equation*}
\begin{split}
	I &= (\pi(a)^3+2\pi(a)^2\pi(b), \pi(a)\pi(b)^2, \pi(a)^2\pi(b)+\pi(a)\pi(b)^2,\pi(a)\pi(b)^2,\pi(b)^3)\\
	&= (\pi(a)^2 + \pi(a)^2\pi(b), \pi(a)\pi(b)^2,\pi(a)\pi(b),\pi(a)\pi(b)^2,\pi(b)^3),
\end{split}
\end{equation*}
where for the second line we used that $\pi(a)+\pi(b)=1$. 

Recall that all elements in $F(|A|,k)$ of length $k$ are constant maps. We are interested in the probability that an 
element of length $1\le \ell<k$ is a constant map when $F(|A|,k)$ acts on right congruences. This is intuitively related 
to the \defn{hitting time} (or waiting time) to constant map. As we will show in Section~\ref{lrc}, there is a lattice 
structure imposed on the set of right congruences with partial order being inclusion. It turns out that we can 
approximate right congruences by \defn{special right congruences} as introduced in 
Section~\ref{section.special right congruences} using certain meets and joins in this lattice. Special right congruences 
in turn are associated to semaphore codes as defined in Section~\ref{section.semaphore codes}, on which it is easy 
to compute the hitting time (see Section~\ref{section.random walks}). The hitting time of the approximation (given 
by a semaphore code) and the right congruence turn out to be the same, and the approximation is finer than the
right congruence. The stationary distributions of the two are simply related by ``lumping''.

\bigskip

Let us now turn our attention to semaphore codes. For a fixed alphabet $A$, which we assume to be a finite non-empty set,
denote by $A^+$ the set of all strings $a_1 \ldots a_\ell$ of length $\ell \ge 1$ over $A$ with multiplication given
by concatenation. Thus $(A^+,A)$ is the free semigroup with generators $A$ (since every semigroup $(S,\cdot)$
generated by a subset $A\subseteq S$ is a surmorphism of $(A^+,A)$ by mapping $a_1\ldots a_\ell \to
a_1 \cdot a_2 \cdot \ldots \cdot a_\ell \in S$). Furthermore, let $A^* = A^+ \cup \{1\}$, so that $A^*$ is $A^+$
with the identity added; it is the free monoid generated by $A$. The semigroup $A^+$ has three orders:
``is a suffix", ``is a prefix", and ``is a factor". In particular, for $u,v \in A^+$
\begin{equation*}
\begin{split}
	& \text{$u$ is a suffix of $v$} \quad \Longleftrightarrow \quad
	\text{$\exists w \in A^*$ such that $w u = v$,}\\
	& \text{$u$ is a prefix of $v$} \quad \Longleftrightarrow \quad
	\text{$\exists w \in A^*$ such that $u w = v$,}\\
	& \text{$u$ is a factor of $v$} \quad \Longleftrightarrow \quad
	\text{$\exists w_1, w_2 \in A^*$ such that $w_1 u w_2 = v$.}
\end{split}
\end{equation*}
A \defn{suffix code} $C$ of $A^+$ (or over $A$) is a subset $C \subseteq A^+$ so that all elements in $C$ are
pairwise incomparable in the suffix order~\cite{BPR:2010}.

A \defn{semaphore code}~\cite{BPR:2010} is a suffix code $S$ over $A$ for which there is a right action in
the following sense:
\begin{equation}
\label{equation.semaphore}
\begin{split}
&\text{If $u \in S \subseteq A^+$ and $a \in A$, then $u a$ has a suffix in $S$ (and hence a unique
suffix of $u a$).} \\
&\text{The right action $u.a$ is the suffix of $u a$ that is in $S$.}
\end{split}
\end{equation}
(The dual concept of prefix codes and left actions is often used in the literature, see for example~\cite{BPR:2010}).
For example, $S=\{ ba^j \mid j \ge 0\} =: ba^*$ is an infinite semaphore code with right action
\[
	ba^j.a = ba^{j+1} \qquad \text{and} \qquad ba^j.b = b.
\]
In practice, to check whether a suffix code is a semaphore code one merely needs to check the first line 
of~\eqref{equation.semaphore}. For example, $C=\{a,bb\}$ is a suffix code, but $a.b$ has no suffix in $C$, so
that $C$ is not a semaphore code.

Semaphore codes over $A$ are inherently related to ideals of $A^+$. A subset $I \subseteq A^+$ is an
\defn{ideal} if $u I v \subseteq I$ for all $u,v \in A^*$. Similarly, $L \subseteq A^+$ is a
\defn{left ideal} if $u L \subseteq L$ for all $u \in A^*$. In this setting, suffix codes over $A$
are precisely the suffix minimal elements of a left ideal $L$.

Now given an ideal $I \subseteq A^+$ we construct a semaphore code as follows. Given
$u = a_j \ldots a_2 a_1 \in A^+$, check whether $u$ is in $I$. If $u \not \in I$, ignore $u$.
If $u \in I$, we find the (necessarily unique) index $1\le i \le j$ such that $a_{i-1} \ldots a_1 \not \in I$,
but $a_i \ldots a_1 \in I$. Then $a_i \ldots a_1$ is a code word and the set of all such words forms the 
semaphore code $S=:I \beta_\ell$, as can be readily verified. It is easy to show that
\[
	I \longleftrightarrow I \beta_\ell
\]
is a bijection between ideals $I \subseteq A^+$ and semaphore codes over $A$, see 
Proposition~\ref{semaphore}. Hence semaphore codes are precisely the suffix minimal elements
of an ideal $I \subseteq A^+$. Since ideals are ubiquitous in mathematics, so are semaphore codes!

\bigskip

As mentioned earlier, the set of right congruences $\rc(A^k)$ is a finite lattice under the inclusion order on
the congruence classes, where the meet is given by intersection. We prove that $\rc(A^k)$ is semimodular, but 
not modular in general, and thus satisfies the Jordan--Dedekind condition that all maximal chains are of the same 
length. Also for $|A|\ge 2$ and $k\ge 2$, $\rc(A^k)$ is not generated by its atoms. See Section~\ref{lrc} for more
details.

Denote by $\sem(A^k)$ the set of semaphore codes coming from ideals $I \supseteq A^k$.
This means that all codewords of $\sem(A^k)$ have length less than or equal to $k$ (so the code
is finite) and every member of $A^k$ has a suffix in the code. Starting with a semaphore code $S$
and restricting the codewords of $S$ to those of length $\le k$, might not yield a finite semaphore code.
But it is always possible to add codewords of length $k$ to this length restricted semaphore code to
obtain $S_k \in \sem(A^k)$. This process of adding codewords of length $k$ which have no suffix 
in the restricted words is unique. For example, we have seen that $S = ba^*$ is a semaphore code.
If we take $k=3$, we obtain $\{b,ba,ba^2\}$. However, $aaa$ has no suffix in this set, so it needs to be
added to obtain the restricted semaphore code $S_3 = \{b,ba,baa,aaa\}$.
In~\cite{RSS:2016} we show that if $S$ is a semaphore code, then the finite
semaphore code $S_k$ converges to $S$ in some precise sense.

Now each semaphore code $S \in \sem(A^k)$ gives a right congruence $\rho \in \rc(A^k)$ as follows:
\begin{equation}
\label{equation.S src}
	\text{For two strings $u,v\in A^k$, we say $u \sim_S v$ if
	$u$ and $v$ have a common suffix in $S$.}
\end{equation}
It is not too hard to verify that $\sim_S$ defines a right congruence on $A^k$. For example, for $A=\{a,b\}$
\[
	S =\{aa,ab,aba,bba,abb,bbb\} \in \sem_3(A)
\]
yields the right congruence in $\rc(A^3)$
\begin{equation}
\label{equation.src example}
	\{aaa,baa\}, \{aab,bab\}, \{aba\}, \{bba\}, \{abb\}, \{bbb\}.
\end{equation}
We denote all elements of $\rc(A^k)$  that arise from semaphore codes in $\sem(A^k)$ by
$\src(A^k)$, the \defn{special right congruences} of $\rc(A^k)$. We prove in Section~\ref{section.special right congruences}
that $\src(A^k)$ is a full (meaning that top and bottom agree) sublattice of $\rc(A^k)$, so that each element
$\rho \in \rc(A^k)$ has a \defn{unique largest lower (finer) approximation} denoted by $\underline{\rho}$, 
namely $\underline{\rho}$ is the join of all elements in $\src(A^k)$ contained in $\rho$.
We will also prove in Section~\ref{section.special right congruences}, and the reader can verify this, 
that the right congruence in~\eqref{equation.rc example} is not a special right congruence, but the special
right congruence in~\eqref{equation.src example} is the unique lower approximation. 
 
\bigskip

As for the de Bruijn graphs, we have random walks on semaphore codes since there is a right action
of a semigroup on semaphore codes. If $S$ is a semaphore code over the alphabet $A$ and
$\pi \colon A \to [0,1]$ is any probability distribution on $A$, namely $\sum_{a\in A} \pi(a)=1$, 
then~\cite[Proposition 3.5.1]{BPR:2010}
\[
	\sum_{s\in S} \pi(s) =1,
\]
where $\pi(s) = \pi(a_1) \cdots \pi(a_\ell)$ if $s=a_1 \ldots a_\ell$. This means in particular that 
$S$ is a \defn{maximal code} with respect to inclusion.

We can now construct a random walk with state space given by the code words in $S$ using the right
action given in~\eqref{equation.semaphore}. Defining the $|S|\times |S|$ monomial matrix $\mathcal{T}(a)$
for each $a\in A$ by $\mathcal{T}(a)_{s,s.a}=1$ and $0$ otherwise for all $s\in S$, we obtain the 
transition matrix as
\[
	\mathcal{T} = \sum_{a\in A} \pi(a) \mathcal{T}(a).
\]
We prove in Theorem~\ref{theorem.stationary} that the stationary distribution $I$ of $\mathcal{T}$ is given
by $I=(\pi(s))_{s\in S}$. Furthermore, the probability that a word of length $\ell$ is a reset (or constant map) is
\[
	P(\ell) = \sum_{\substack{s\in S\\ \ell(s) \le \ell}} \pi(s),
\]
see Theorem~\ref{theorem.length prob}. This probability is related to the hitting time to reset. 
For example, for the semaphore code $S = ba^*$, all words $w$ are resets unless $w=a^\ell$.
The probability that a string of length 3 is a reset is $P(3) = \pi(b) + \pi(b) \pi(a) + \pi(b) \pi(a)^2
=1-\pi(a)^3$. For more details see Section~\ref{section.random walks}.

\bigskip

We are now able to give a more direct construction of the special right congruence $\underline{\rho}$
for $\rho \in \rc(A^k)$, the best lower approximation of $\rho$ in $\src(A^k)$. Define
\[
	\res(\rho) = \{ w \in A^+ \mid \text{$w$ is a reset on $A^k/\rho$} \}.
\]
Then we prove that $\res(\rho)$ is an ideal of $A^k \subseteq A^+$ and the special right congruence
associated to the semaphore code given by this ideal is $\underline{\rho}$. An immediate consequence is
that $\rho$ and $\underline{\rho}$ have the same hitting time to reset, but in general different stationary
distributions. In general, $\underline{\rho}$ has more congruence classes than $\rho$, so the stationary
distributions cannot be the same. Note that both distributions are determined by lumping from the 
product distribution of the de Bruijn random walk on $A^k$.
In applications a metric is placed on all distributions of $\rc(A^k)$. Then the probability
distribution $\pi$ on $A$ is chosen such that the distance between $I_\rho$ and $I_{\underline{\rho}}$
is minimal. This is called the principle of choosing a ``correct" or ``good" probability distribution $\pi$ on
$A$. 

\bigskip

The paper is organized as follows. In Section~\ref{section.algebra} we provide the algebraic background
of the semigroups related to right congruences. The precise definition of resets is given in
Section~\ref{regr}. Semaphore codes are introduced in Section~\ref{section.semaphore codes}.
In Sections~\ref{rcmi} right congruence and their properties are studied, in particular the lattice structure
in Section~\ref{lrc}. Special right congruences are the subject of Section~\ref{section.special right congruences}.
Random walks on semaphore codes are studied in Section~\ref{section.random walks}.
Note that the semaphore codes introduced in Section~\ref{section.semaphore codes} can be infinite.
The analysis in terms of random walks in Section~\ref{section.random walks} is valid for both
finite and infinite semaphore codes. In all other sections we restrict to finite lengths words and
codes.

\subsection*{Acknowledgements}

We would like to thank Arvind Ayyer, Benjamin Steinberg and Nicolas M. Thi\'ery for discussions.

The first author thanks the Simons Foundation--Collaboration Grants for Mathematicians for  travel grant $\#313548$.
The second author was partially supported by NSF grants OCI--1147247 and DMS--1500050.
The third author was partially supported by CMUP (UID/MAT/00144/2013),
which is funded by FCT (Portugal) with national (MEC) and European
structural funds through the programs FEDER, under the partnership agreement PT2020.

\section{Algebraic foundations}
\label{section.algebra}

\subsection{Elliptic maps on rooted trees}
\label{subsection.elliptic}

Elliptic maps on finite trees were considered by Rhodes and Silva~\cite{Rhodes:1991, RS:2012}.
A \defn{tree} is a connected graph that does not contain a closed walk in which all vertices are distinct.
A \defn{leaf} of a tree is a vertex of degree 1, that is, a vertex that connects to exactly one edge.
A \defn{rooted tree} is a tree in which a particular node is designated as the root. In this case,
if a vertex $u$ is on the path from the root to another vertex $v$, we say that $u$ is an \defn{ancestor} of $v$,
or equivalently, that $v$ is a \defn{descendant} of $u$. If $u$ and $v$ are
adjacent, we say that $u$ is the \defn{parent} of $v$, which is the \defn{child} of $u$.

Given a rooted tree $T$, we denote by $\Vert(T)$ the set of vertices
of $T$. The distance between two vertices is the minimum number of
edges in a path between them. An \defn{elliptic map} on $T$ is a
mapping $\Vert(T) \to \Vert(T)$ preserving adjacency and distance to
the root. Equivalently, an elliptic map on $T$ is a contraction
(decreases distances between vertices) while preserving distance to the
root, or a mapping fixing the root and preserving parenthood.
We shall write functions on the right since we will deal with
right actions and compositions.
Elliptic maps on a fixed rooted tree form a monoid under composition.

\begin{figure}[t]
\begin{center}
{ \newcommand{\nodea}{\node[draw,circle] (a) {$$}
;}\newcommand{\nodeb}{\node[draw,circle] (b) {$$}
;}\newcommand{\nodec}{\node[draw,circle] (c) {$$}
;}\newcommand{\noded}{\node[draw,circle] (d) {$$}
;}\newcommand{\nodee}{\node[draw,circle] (e) {$$}
;}\newcommand{\nodef}{\node[draw,circle] (f) {$$}
;}\newcommand{\nodeg}{\node[draw,circle] (g) {$$}
;}\newcommand{\nodeh}{\node[draw,circle] (h) {$$}
;}\newcommand{\nodei}{\node[draw,circle] (i) {$$}
;}\begin{tikzpicture}[auto]
\matrix[column sep=.3cm, row sep=.3cm,ampersand replacement=\&]{
         \&         \&         \& \nodea  \&         \&         \&         \\ 
         \& \nodeb  \&         \&         \&         \& \nodef  \&         \\ 
 \nodec  \& \noded  \& \nodee  \&         \& \nodeg  \& \nodeh  \& \nodei  \\
};

\path[ultra thick, red] (b) edge (c) edge (d) edge (e)
	(f) edge (g) edge (h) edge (i)
	(a) edge (b) edge (f);
\end{tikzpicture}}
\end{center}
\caption{Rooted tree $T(2,3)$.
\label{figure.rooted tree}}
\end{figure}
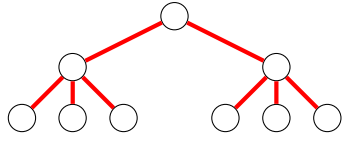

\begin{figure}[t]
\begin{center}
{ \newcommand{\nodea}{\node[draw,circle] (a) {$$}
;}\newcommand{\nodeb}{\node[draw,circle] (b) {$$}
;}\newcommand{\nodec}{\node[draw,circle] (c) {$$}
;}\newcommand{\noded}{\node[draw,circle] (d) {$$}
;}\newcommand{\nodee}{\node[draw,circle] (e) {$$}
;}\newcommand{\nodef}{\node[draw,circle] (f) {$$}
;}\newcommand{\nodeg}{\node[draw,circle] (g) {$$}
;}\newcommand{\nodeh}{\node[draw,circle] (h) {$$}
;}\newcommand{\nodei}{\node[draw,circle] (i) {$$}
;}\begin{tikzpicture}[auto]
\matrix[column sep=.3cm, row sep=.3cm,ampersand replacement=\&]{
         \&         \&         \& \nodea  \&         \&         \&         \\ 
         \& \nodeb  \&         \&         \&         \& \nodef  \&         \\ 
 \nodec  \& \noded  \& \nodee  \&         \& \nodeg  \& \nodeh  \& \nodei  \\
};

\node at (-0.5,0.7) {$r_0$};
\node at (-1.8,0) {$v_1$};
\node at (1.8,0) {$v_2$};
\node at (-2.2,-1.1) {$v_{11}$};
\node at (-1.4,-1.1) {$v_{12}$};
\node at (-0.55,-1.1) {$v_{13}$};
\node at (0.55,-1.1) {$v_{21}$};
\node at (1.4,-1.1) {$v_{22}$};
\node at (2.2,-1.1) {$v_{23}$};

\path[ultra thick, red] (b) edge (c) edge (d) edge (e)
	(f) edge (g) edge (h) edge (i)
	(a) edge (b) edge (f);
\end{tikzpicture}}
\raisebox{1cm}{$\qquad \mapsto \qquad$}
\raisebox{-0.2cm}{
{ \newcommand{\nodea}{\node[draw,circle] (a) {$$}
;}\newcommand{\nodeb}{\node[draw,circle] (b) {$$}
;}\newcommand{\nodec}{\node[draw,circle] (c) {$$}
;}\newcommand{\noded}{\node[draw,circle] (d) {$$}
;}\newcommand{\nodee}{\node[draw,circle] (e) {$$}
;}\newcommand{\nodef}{\node[draw,circle] (f) {$$}
;}\newcommand{\nodeg}{\node[draw,circle] (g) {$$}
;}\newcommand{\nodeh}{\node[draw,circle] (h) {$$}
;}\newcommand{\nodei}{\node[draw,circle] (i) {$$}
;}\begin{tikzpicture}[auto]
\matrix[column sep=.3cm, row sep=.3cm,ampersand replacement=\&]{
         \&         \&         \& \nodea  \&         \&         \&         \\ 
         \& \nodeb  \&         \&         \&         \& \nodef  \&         \\ 
 \nodec  \& \noded  \& \nodee  \&         \& \nodeg  \& \nodeh  \& \nodei  \\
};

\node at (-0.5,0.7) {$r_0$};
\node at (-1.8,0) {$v_2$};
\node at (1.8,0) {$v_1$};
\node at (-2.2,-1.1) {$v_{22}$};
\node at (-2.2,-1.5) {$v_{23}$};
\node at (-1.4,-1.1) {$v_{21}$};
\node at (0.55,-1.1) {$v_{11}$};
\node at (0.55,-1.5) {$v_{12}$};
\node at (2.2,-1.1) {$v_{13}$};

\path[ultra thick, red] (b) edge (c) edge (d) edge (e)
	(f) edge (g) edge (h) edge (i)
	(a) edge (b) edge (f);
\end{tikzpicture}}}
\end{center}
\caption{Elliptic map $\varphi \colon \Vert(T) \to \Vert(T)$ on $T:=T(2,3)$ which maps $r_0\mapsto r_0$,
$v_1 \mapsto v_2$, $v_2\mapsto v_1$, $v_{11} \mapsto v_{21}$, $v_{12} \mapsto v_{21}$, $v_{13} \mapsto v_{23}$,
$v_{21} \mapsto v_{12}$, $v_{22} \mapsto v_{11}$, $v_{23} \mapsto v_{11}$.
\label{figure.elliptic map}}
\end{figure}
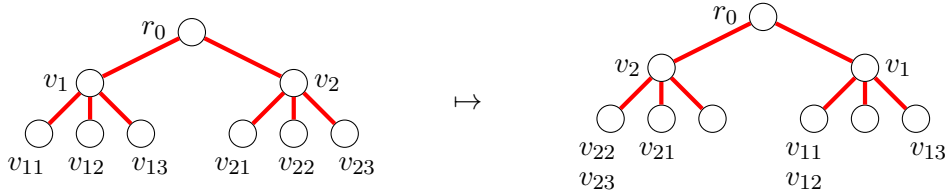

Let $T:=T(n_0,\ldots,n_N)$ be a uniformly branching rooted tree, where
all leaves are at distance $N+1$ from  
the root $r_0$ and each vertex at distance (or level) $k$ from the
root has $n_k$ children for $k = 0,\ldots,N$. 
An example of a uniformly branching rooted tree is given in
Figure~\ref{figure.rooted tree}. 
An example of an elliptic map on this tree is given in Figure~\ref{figure.elliptic map}.

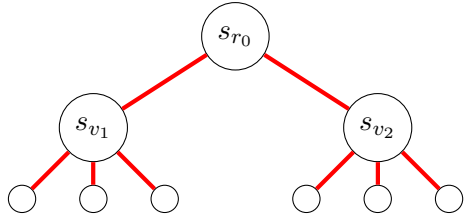
\begin{figure}[t]
\begin{center}
{ \newcommand{\nodea}{\node[draw,circle] (a) {$s_{r_0}$}
;}\newcommand{\nodeb}{\node[draw,circle] (b) {$s_{v_1}$}
;}\newcommand{\nodec}{\node[draw,circle] (c) {$$}
;}\newcommand{\noded}{\node[draw,circle] (d) {$$}
;}\newcommand{\nodee}{\node[draw,circle] (e) {$$}
;}\newcommand{\nodef}{\node[draw,circle] (f) {$s_{v_2}$}
;}\newcommand{\nodeg}{\node[draw,circle] (g) {$$}
;}\newcommand{\nodeh}{\node[draw,circle] (h) {$$}
;}\newcommand{\nodei}{\node[draw,circle] (i) {$$}
;}\begin{tikzpicture}[auto]
\matrix[column sep=.3cm, row sep=.3cm,ampersand replacement=\&]{
         \&         \&         \& \nodea  \&         \&         \&         \\ 
         \& \nodeb  \&         \&         \&         \& \nodef  \&         \\ 
 \nodec  \& \noded  \& \nodee  \&         \& \nodeg  \& \nodeh  \& \nodei  \\
};

\path[ultra thick, red] (b) edge (c) edge (d) edge (e)
	(f) edge (g) edge (h) edge (i)
	(a) edge (b) edge (f);
\end{tikzpicture}}
\end{center}
\caption{Elliptic map of Figure~\ref{figure.elliptic map} written with component actions: $s_{r_0}$ is the map
$v_1 \mapsto v_2$, $v_2 \mapsto v_1$, $s_{v_1}$ is the map $v_{11} \mapsto v_{21}$, $v_{12} \mapsto v_{21}$,
$v_{13} \mapsto v_{23}$, and $s_{v_2}$ is the map $v_{21} \mapsto v_{12}$, $v_{22} \mapsto v_{11}$,
$v_{23} \mapsto v_{11}$.
\label{figure.component action}}
\end{figure}

\begin{figure}[t]
\begin{center}
{ \newcommand{\nodea}{\node[draw,circle] (a) {$$}
;}\newcommand{\nodeb}{\node[draw,circle] (b) {$$}
;}\newcommand{\nodec}{\node[draw,circle] (c) {$$}
;}\newcommand{\noded}{\node[draw,circle] (d) {$$}
;}\newcommand{\nodee}{\node[draw,circle] (e) {$$}
;}\newcommand{\nodef}{\node[draw,circle] (f) {$$}
;}\newcommand{\nodeg}{\node[draw,circle] (g) {$$}
;}\newcommand{\nodeh}{\node[draw,circle] (h) {$$}
;}\newcommand{\nodei}{\node[draw,circle] (i) {$$}
;}\newcommand{\nodej}{\node[draw,circle] (j) {$$}
;}\newcommand{\nodeba}{\node[draw,circle] (ba) {$$}
;}\newcommand{\nodebb}{\node[draw,circle] (bb) {$$}
;}\newcommand{\nodebc}{\node[draw,circle] (bc) {$$}
;}\newcommand{\nodebd}{\node[draw,circle] (bd) {$$}
;}\newcommand{\nodebe}{\node[draw,circle] (be) {$$}
;}\newcommand{\nodebf}{\node[draw,circle] (bf) {$s_w$}
;}\newcommand{\nodebg}{\node[draw,circle] (bg) {$$}
;}\newcommand{\nodebh}{\node[draw,circle] (bh) {$$}
;}\newcommand{\nodebi}{\node[draw,circle] (bi) {$$}
;}\newcommand{\nodebj}{\node[draw,circle] (bj) {$$}
;}\newcommand{\nodeca}{\node[draw,circle] (ca) {$$}
;}\newcommand{\nodecb}{\node[draw,circle] (cb) {$$}
;}\newcommand{\nodecc}{\node[draw,circle] (cc) {$$}
;}\newcommand{\nodecd}{\node[draw,circle] (cd) {$$}
;}\newcommand{\nodece}{\node[draw,circle] (ce) {$$}
;}\newcommand{\nodecf}{\node[draw,circle] (cf) {$$}
;}\newcommand{\nodecg}{\node[draw,circle] (cg) {$$}
;}\begin{tikzpicture}[auto]
\matrix[column sep=.3cm, row sep=.3cm,ampersand replacement=\&]{
         \&         \&         \&         \&         \&         \&         \&         \&         \& \nodea  \&         \&         \&         \&         \&         \&         \&         \&         \&         \\ 
         \&         \&         \&         \& \nodeb  \&         \&         \&         \&         \&         \&         \&         \&         \&         \& \nodebe \&         \&         \&         \&         \\ 
         \& \nodec  \&         \&         \& \nodeg  \&         \&         \& \nodeba \&         \&         \&         \& \nodebf \&         \&         \& \nodebj \&         \&         \& \nodecd \&         \\ 
 \noded  \& \nodee  \& \nodef  \& \nodeh  \& \nodei  \& \nodej  \& \nodebb \& \nodebc \& \nodebd \&         \& \nodebg \& \nodebh \& \nodebi \& \nodeca \& \nodecb \& \nodecc \& \nodece \& \nodecf \& \nodecg \\
};
\node at (0.7,-0.3) {$w$};
\node at (0.5,-1.6) {$v_1$};
\node at (1.4,-1.6) {$v_2$};
\node at (2.3,-1.6) {$v_3$};

\path[ultra thick, red] 
        (c) edge (d) edge (e) edge (f)
	(g) edge (h) edge (i) edge (j)
	(ba) edge (bb) edge (bc) edge (bd)
	(b) edge (c) edge (g) edge (ba)
	(bf) edge (bg) edge (bh) edge (bi)
	(bj) edge (ca) edge (cb) edge (cc)
	(cd) edge (ce) edge (cf) edge (cg)
	(be) edge (bf) edge (bj) edge (cd)
	(a) edge (b) edge (be);
\end{tikzpicture}}
\end{center}
\caption{Component action at vertex $w$ of an elliptic map on $T(2,3,3)$. The component action $s_w$ is a map on the 
children of $w$, namely on $\{v_1,v_2,v_3\}$, and maps into the children of the image of $w$ under the elliptic
map.
\label{figure.m}}
\end{figure}
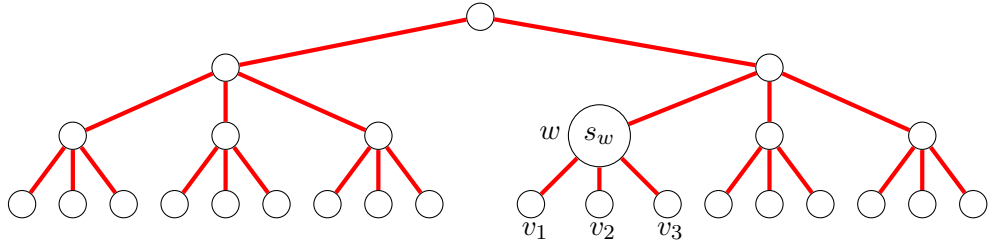

There is another way to represent an elliptic map $\varphi$ using \defn{component actions}. 
Namely, a given vertex $v\in \Vert(T)$ at level $k$ is completely specified by the unique path 
$r_0 \to w_1 \to \cdots \to w_k=v$ from the root. Since elliptic maps
preserve parenthood, the image of this path
under the elliptic map $r_0 \to (w_1)\varphi \to \cdots \to (w_k)\varphi = (v)\varphi$ is again a path,
this time from $r_0$ to $(v)\varphi$. Hence $\varphi$ can be defined recursively:
given the map from path $r_0 \to w_1 \to \cdots \to w_{k-1}$ to
$r_0 \to (w_1)\varphi \to \cdots \to (w_{k-1})\varphi$, we can define a map $s_w$ from the children of
$w:=w_{k-1}$ to the children of $(w_{k-1})\varphi$. The map $s_w$ is called the component action at
vertex $w$. Graphically, we place $s_w$ on the vertex $w$ for every vertex $w$ that is not a leaf.
See Figure~\ref{figure.m}. The elliptic map of Figure~\ref{figure.elliptic map} is written using component
actions in Figure~\ref{figure.component action}.

\begin{figure}[t]
\begin{center}
{ \newcommand{\nodea}{\node[draw,circle] (a) {$s_{r_0}$}
;}\newcommand{\nodeb}{\node[draw,circle] (b) {$s_{v_1}$}
;}\newcommand{\nodec}{\node[draw,circle] (c) {$$}
;}\newcommand{\noded}{\node[draw,circle] (d) {$$}
;}\newcommand{\nodee}{\node[draw,circle] (e) {$$}
;}\newcommand{\nodef}{\node[draw,circle] (f) {$s_{v_2}$}
;}\newcommand{\nodeg}{\node[draw,circle] (g) {$$}
;}\newcommand{\nodeh}{\node[draw,circle] (h) {$$}
;}\newcommand{\nodei}{\node[draw,circle] (i) {$$}
;}\begin{tikzpicture}[auto]
\matrix[column sep=.3cm, row sep=.3cm,ampersand replacement=\&]{
         \&         \&         \& \nodea  \&         \&         \&         \\ 
         \& \nodeb  \&         \&         \&         \& \nodef  \&         \\ 
 \nodec  \& \noded  \& \nodee  \&         \& \nodeg  \& \nodeh  \& \nodei  \\
};

\path[ultra thick, red] (b) edge (c) edge (d) edge (e)
	(f) edge (g) edge (h) edge (i)
	(a) edge (b) edge (f);
\end{tikzpicture}}
\raisebox{1.5cm}{$\circ$}
{ \newcommand{\nodea}{\node[draw,circle] (a) {$t_{r_0}$}
;}\newcommand{\nodeb}{\node[draw,circle] (b) {$t_{v_1}$}
;}\newcommand{\nodec}{\node[draw,circle] (c) {$$}
;}\newcommand{\noded}{\node[draw,circle] (d) {$$}
;}\newcommand{\nodee}{\node[draw,circle] (e) {$$}
;}\newcommand{\nodef}{\node[draw,circle] (f) {$t_{v_2}$}
;}\newcommand{\nodeg}{\node[draw,circle] (g) {$$}
;}\newcommand{\nodeh}{\node[draw,circle] (h) {$$}
;}\newcommand{\nodei}{\node[draw,circle] (i) {$$}
;}\begin{tikzpicture}[auto]
\matrix[column sep=.3cm, row sep=.3cm,ampersand replacement=\&]{
         \&         \&         \& \nodea  \&         \&         \&         \\ 
         \& \nodeb  \&         \&         \&         \& \nodef  \&         \\ 
 \nodec  \& \noded  \& \nodee  \&         \& \nodeg  \& \nodeh  \& \nodei  \\
};

\path[ultra thick, red] (b) edge (c) edge (d) edge (e)
	(f) edge (g) edge (h) edge (i)
	(a) edge (b) edge (f);
\end{tikzpicture}}
\raisebox{1.5cm}{$=$}
{ \newcommand{\nodea}{\node[draw,circle] (a) {$s_{r_0} t_{r_0}$}
;}\newcommand{\nodeb}{\node[draw,circle] (b) {$s_{v_1} t_{(v_1)s_{r_0}}$}
;}\newcommand{\nodec}{\node[draw,circle] (c) {$$}
;}\newcommand{\noded}{\node[draw,circle] (d) {$$}
;}\newcommand{\nodee}{\node[draw,circle] (e) {$$}
;}\newcommand{\nodef}{\node[draw,circle] (f) {$s_{v_2} t_{(v_2)s_{r_0}}$}
;}\newcommand{\nodeg}{\node[draw,circle] (g) {$$}
;}\newcommand{\nodeh}{\node[draw,circle] (h) {$$}
;}\newcommand{\nodei}{\node[draw,circle] (i) {$$}
;}\begin{tikzpicture}[auto]
\matrix[column sep=.3cm, row sep=.3cm,ampersand replacement=\&]{
         \&         \&         \& \nodea  \&         \&         \&         \\ 
         \& \nodeb  \&         \&         \&         \& \nodef  \&         \\ 
 \nodec  \& \noded  \& \nodee  \&         \& \nodeg  \& \nodeh  \& \nodei  \\
};

\path[ultra thick, red] (b) edge (c) edge (d) edge (e)
	(f) edge (g) edge (h) edge (i)
	(a) edge (b) edge (f);
\end{tikzpicture}}
\end{center}
\caption{Composition or product of two elliptic maps on the rooted tree in Figure~\ref{figure.rooted tree}.
\label{figure.product}}
\end{figure}
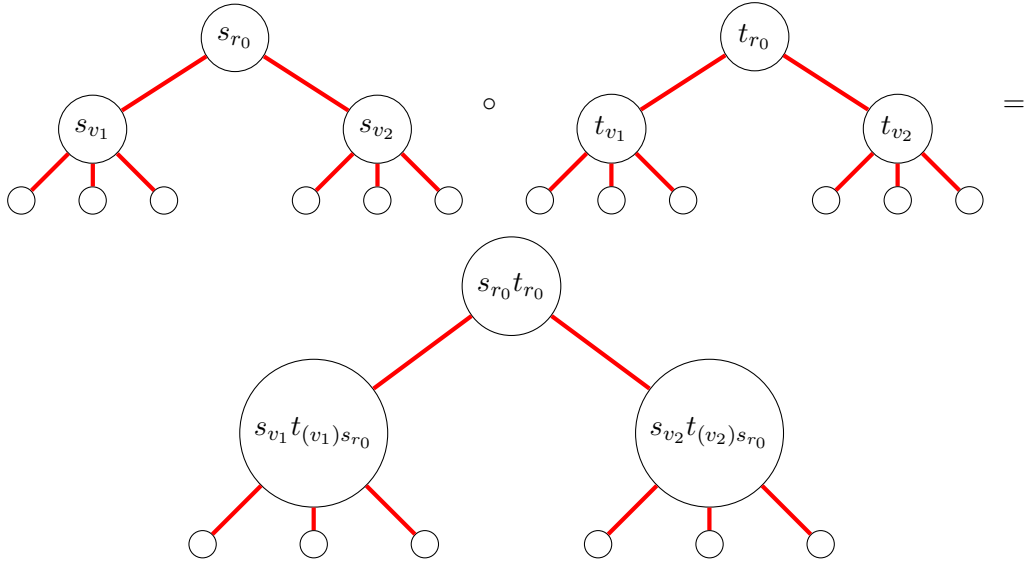

As mentioned before, the \defn{product of elliptic maps} is composition, which is another elliptic map. 
We can formulate this in terms of the component actions. Let $\varphi$ and $\psi$ be elliptic maps on the same rooted 
tree $T$ with component action $s_v$ and $t_v$ at vertex $v\in \Vert(T)$ that is not a leaf, respectively. Then the 
component action of $\varphi \circ \psi$ at vertex $v$ is $s_v t_{(v)s_w}$, where $w$ is the parent of $v$.
An example is given in Figure~\ref{figure.product}.

\medskip

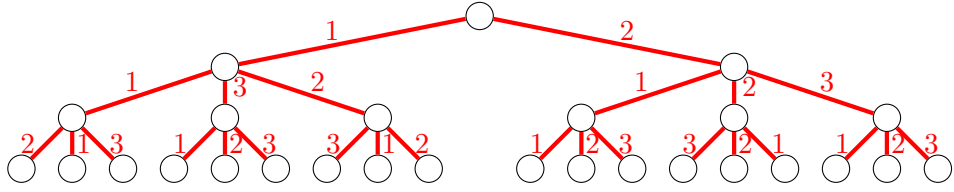
\begin{figure}[t]
\begin{center}
{ \newcommand{\nodea}{\node[draw,circle] (a) {$$}
;}\newcommand{\nodeb}{\node[draw,circle] (b) {$$}
;}\newcommand{\nodec}{\node[draw,circle] (c) {$$}
;}\newcommand{\noded}{\node[draw,circle] (d) {$$}
;}\newcommand{\nodee}{\node[draw,circle] (e) {$$}
;}\newcommand{\nodef}{\node[draw,circle] (f) {$$}
;}\newcommand{\nodeg}{\node[draw,circle] (g) {$$}
;}\newcommand{\nodeh}{\node[draw,circle] (h) {$$}
;}\newcommand{\nodei}{\node[draw,circle] (i) {$$}
;}\newcommand{\nodej}{\node[draw,circle] (j) {$$}
;}\newcommand{\nodeba}{\node[draw,circle] (ba) {$$}
;}\newcommand{\nodebb}{\node[draw,circle] (bb) {$$}
;}\newcommand{\nodebc}{\node[draw,circle] (bc) {$$}
;}\newcommand{\nodebd}{\node[draw,circle] (bd) {$$}
;}\newcommand{\nodebe}{\node[draw,circle] (be) {$$}
;}\newcommand{\nodebf}{\node[draw,circle] (bf) {$$}
;}\newcommand{\nodebg}{\node[draw,circle] (bg) {$$}
;}\newcommand{\nodebh}{\node[draw,circle] (bh) {$$}
;}\newcommand{\nodebi}{\node[draw,circle] (bi) {$$}
;}\newcommand{\nodebj}{\node[draw,circle] (bj) {$$}
;}\newcommand{\nodeca}{\node[draw,circle] (ca) {$$}
;}\newcommand{\nodecb}{\node[draw,circle] (cb) {$$}
;}\newcommand{\nodecc}{\node[draw,circle] (cc) {$$}
;}\newcommand{\nodecd}{\node[draw,circle] (cd) {$$}
;}\newcommand{\nodece}{\node[draw,circle] (ce) {$$}
;}\newcommand{\nodecf}{\node[draw,circle] (cf) {$$}
;}\newcommand{\nodecg}{\node[draw,circle] (cg) {$$}
;}\begin{tikzpicture}[auto]
\matrix[column sep=.3cm, row sep=.3cm,ampersand replacement=\&]{
         \&         \&         \&         \&         \&         \&         \&         \&         \& \nodea  \&         \&         \&         \&         \&         \&         \&         \&         \&         \\ 
         \&         \&         \&         \& \nodeb  \&         \&         \&         \&         \&         \&         \&         \&         \&         \& \nodebe \&         \&         \&         \&         \\ 
         \& \nodec  \&         \&         \& \nodeg  \&         \&         \& \nodeba \&         \&         \&         \& \nodebf \&         \&         \& \nodebj \&         \&         \& \nodecd \&         \\ 
 \noded  \& \nodee  \& \nodef  \& \nodeh  \& \nodei  \& \nodej  \& \nodebb \& \nodebc \& \nodebd \&         \& \nodebg \& \nodebh \& \nodebi \& \nodeca \& \nodecb \& \nodecc \& \nodece \& \nodecf \& \nodecg \\
};

\path[ultra thick, red] 
        (c) edge node [left]  {$2$} (d) edge node [right]  {\hspace{-.1cm}$1$} (e) edge node [right]  {$3$} (f)
	(g) edge node [left]  {$1$} (h) edge node [right]  {\hspace{-.1cm}$2$} (i) edge node [right]  {$3$} (j)
	(ba) edge node [left]  {$3$} (bb) edge node [right]  {\hspace{-.1cm}$1$} (bc) edge node [right]  {$2$} (bd)
	(b) edge node [near start, left]  {$1\quad$} (c) edge node [near start, right]  {\hspace{-.05cm}$3$} (g) edge node [near start, right]  {$\quad 2$} (ba)
	(bf) edge node [left]  {$1$} (bg) edge node [right]  {\hspace{-.1cm}$2$} (bh) edge node [right]  {$3$} (bi)
	(bj) edge node [left]  {$3$} (ca) edge node [right]  {\hspace{-.1cm}$2$} (cb) edge node [right]  {$1$} (cc)
	(cd) edge node [left]  {$1$} (ce) edge node [right]  {\hspace{-.1cm}$2$} (cf) edge node [right]  {$3$} (cg)
	(be) edge node [near start, left]  {$1\quad$} (bf) edge node [near start, right]  {\hspace{-.05cm}$2$} (bj) edge node [near start, right]  {$\quad 3$} (cd)
	(a) edge node [near start, left]  {$1\qquad$} (b) edge node [near start, right]  {$\qquad 2$} (be);
\end{tikzpicture}}
\end{center}
\caption{Labelled rooted tree $T(2,3,3)$ with labeling sets $X_0=\{1,2\}$, $X_1= X_2= \{1,2,3\}$.
\label{figure.X}}
\end{figure}

Note that a child $v$ of a vertex $w$ can be uniquely specified by the edge $e$ that leads to it.
Hence the path $r_0=w_0 \to w_1 \to \cdots \to w_k=v$ from $r_0$ to $v$ can alternatively be encoded by a sequence
$e_0 \to e_1 \to \cdots \to e_{k-1}$ of edges, where $e_i$ is the edge from vertex $w_i$ to $w_{i+1}$. For us, it will 
be convenient to keep track of the edges by labelling the $n_\ell$ edges leaving a given vertex at level 
$0\le \ell \le N$ bijectively with elements from a set $X_\ell$ with
$|X_\ell|=n_\ell$. The result is a \defn{labelled rooted tree}. See Figure~\ref{figure.X} for an example.
Note that there are lots of ways to label a rooted tree. Labelling the rooted tree is equivalent to specifying 
a coordinate system. Once the labelling $L$ of $T$ is fixed, a sequence $e_0 \to e_1 \to \cdots \to e_{k-1}$ of edges
is determined by an element $(x_0, x_1, \ldots, x_{k-1}) \in X_0 \times X_1 \times \cdots \times X_{k-1}$.

Given a rooted tree $T(n_0,\ldots,n_N)$ with labels in $X=X_0 \times \cdots \times X_N$, 
elliptic maps can now be expressed using the labels giving rise to the \defn{wreath product}.
The component action at level $k$ is described by a semigroup $S_k$ acting faithfully on the right on $X_k$, 
denoted $(X_k,S_k)$. Then the wreath product $(X_0,S_0)\circ \cdots \circ (X_N,S_N)$ is $(X,S)$, where $S$
is the semigroup with component action at level $k$ in $(X_k,S_k)$. More precisely, $\Pi = (\Pi_0,\ldots,\Pi_N)\in S$ 
if $\Pi_0 \in S_0$, $\Pi_1 \colon X_0 \to S_1$, and generally $\Pi_k \colon X_0 \times \cdots \times X_{k-1} \to S_k$
for $1\le k\le N$, so that for $(x_0,\ldots, x_N) \in X$
\begin{equation}
\label{equation.wreath}
	(x_0, \ldots ,x_N) \Pi = \Bigl(x_0.\Pi_0, x_1.(x_0)\Pi_1,
	x_2. (x_0,x_1)\Pi_2, \ldots, x_N. (x_0,\ldots,x_{N-1})\Pi_N \Bigr) \;.
\end{equation}
The semigroup element $m:= (x_0,\ldots, x_{k-1}) \Pi_k \in S_k$ is the \defn{component action} in the vertex
(or component) specified by $(x_0,\ldots, x_{k-1})$. 

\begin{Remark}
The above arguments show that elliptic maps on uniformly branching trees and wreath products are the same
thing (confirming~\cite[Proposition 3.3]{RS:2012}).
\end{Remark}

\defn{Multiplication of wreath products} is given by composition of the component action~\eqref{equation.wreath}.
Graphically on the level of labelled trees directly, the product $\Pi^g \cdot \Pi^f$ for $\Pi^g, \Pi^f \in (X,S)$
translates to the following:
\begin{enumerate}
\item To determine the value of $\Pi^g \cdot \Pi^f$ at vertex $x = (x_0,\ldots,x_{k-1})$ in the labelled rooted tree, 
go to the corresponding vertex in the tree for $\Pi^g$, keep track of all values at the vertices on the way
and act with the corresponding elements on the vertex vector:
\[
	x^g = \Bigl(x_0.\Pi^g_0, x_1.(x_0)\Pi^g_1, x_2.(x_0,x_1)\Pi^g_2, \ldots, x_k.(x_0,\ldots,x_{k-1})\Pi_k^g \Bigr) \;.
\]
\item Then the entry in vertex $(x_0,\ldots,x_{k-1})$ of $\Pi^g \cdot \Pi^f$
is $(x_0,\ldots,x_{k-1})\Pi_k^g  (x^g_0,\ldots,x^g_{k-1})\Pi_k^f$.
\end{enumerate}

\medskip

One of the main questions is ``how restrained can the component action be''?
See the first half of~\cite{Rhodes:2010} and the introduction to~\cite{RS:2009}.

\medskip

The \defn{Prime Decomposition Theorem} of Krohn and Rhodes~\cite{KR:1965} (see also~\cite{Rhodes:2010}
and~\cite[Chapter 4]{RS:2009}) states that every finite semigroup 
divides an iterated wreath product of its finite simple group divisors and copies of the three element aperiodic 
monoid $U_2$ consisting of two right zeroes and an identity. More precisely, a semigroup $S_1$ divides 
semigroup $S_2$, written $S_1 | S_2$, if $S_1$ is a homomorphic image of a subsemigroup of $S_2$.
In addition, $U_2=\{1,a,b\}$ where $xa=a, xb=b$, and $1x=x1=x$ for all $x \in U_2$. A finite semigroup is aperiodic 
if all of its subgroups are trivial.
Alternatively, the Prime Decomposition Theorem says that the basic building blocks of finite semigroups are the 
finite simple groups and semigroups of constant maps with an adjoined identity.

We say that $I \subseteq S$ is an \defn{ideal} of the semigroup $S$ if
$SI \cup IS \subseteq I$. We write then $I \unlhd S$.
The \defn{kernel} of a semigroup $S$, denoted $\ker(S)$, is the unique
minimal nonempty ideal of $S$. If $S$ is a monoid, its
group of units is the subgroup formed by all the invertible
elements. Both kernel and group of units play a major role in this
context.

Let $S_1$ and $S_2$ be semigroups and let $\varphi$ be a homomorphism of $S_1$ into endomorphisms of $S_2$. 
Then the semigroup $S_1 \times_\varphi S_2$ is the \defn{semidirect product} of $S_1$ by $S_2$ with connecting 
homomorphism $\varphi$ (see also~\cite[Section 1.2.2, pg. 23]{RS:2009}). More precisely, $S_1 \times_\varphi S_2$ 
has elements in $S_1 \times S_2$ with multiplication given by 
\[ 
	(s_1,s_2)\cdot (s_1',s_2') =  (s_1 s_1', s_2((s_1')\varphi) \; s_2')\;.
\]
Notice that wreath products are a special case of semidirect products.  In fact, wreath products are
``generic'' semidirect products. Namely up to pseudovarieties, semidirect products, wreath products,
and elliptic products yield the same thing. See~\cite{RS:2009} for all details.

A semigroup $S$ is called \defn{irreducible} if for all finite semigroups $S_1$ and $S_2$ and all connecting 
homomorphisms $\varphi$, $S\mid S_1 \times_\varphi S_2$ implies $S | S_1$ or $S | S_2$. Krohn and
Rhodes~\cite{KR:1965} showed that $S$ is irreducible if and only if either
(a) $S$ is a nontrivial simple group; or
(b) $S$ is one of the four divisors of $U_2$.

\medskip

A \defn{pseudovariety} is a collection of finite semigroups closed under taking finite direct products 
and divisors (that is, subsemigroups and quotients)~\cite{RS:2009}. The monoid $U_2$ is in the 
pseudovariety $\RZ^1$, where $\RZ=[[xy=y]]$ is the pseudovariety of \defn{right zeroes}, meaning that 
all elements $x,y$ in $S\in \RZ$ satisfy the identity $xy=y$. In other
words, $\RZ$ is the pseudovariety generated by semigroups of
constant maps. We denote by $\RZ^1$ the pseudovariety generated by
semigroups of transformations consisting of constant maps plus the
identity mapping. The elements in $\RZ^1$ are also 
called left regular bands, indeed
$\RZ^1 = [[x^2 = x, \; xyx = yx]]$ (cf. \cite[Proposition 7.3.2]{RS:2009}).
Random walks on left regular band are an important new topic~\cite{B:2000, BD:1998}.
This has recently also been generalized to random walks on $\mathscr{R}$-trivial
monoids~\cite{AKS:2014,ASST:2015}.

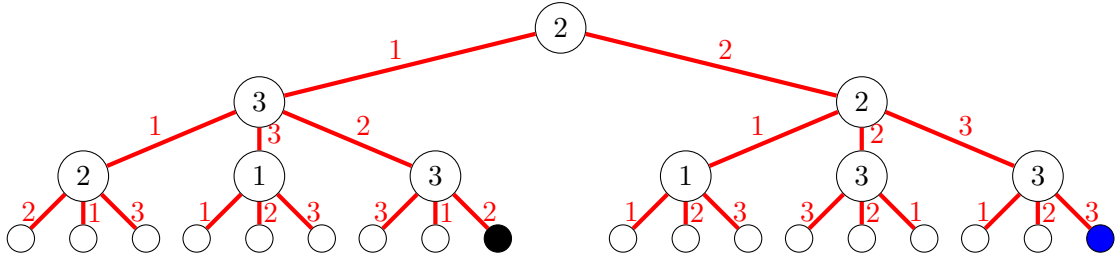
\begin{figure}[t]
\begin{center}
{ \newcommand{\nodea}{\node[draw,circle] (a) {$2$}
;}\newcommand{\nodeb}{\node[draw,circle] (b) {$3$}
;}\newcommand{\nodec}{\node[draw,circle] (c) {$2$}
;}\newcommand{\noded}{\node[draw,circle] (d) {$$}
;}\newcommand{\nodee}{\node[draw,circle] (e) {$$}
;}\newcommand{\nodef}{\node[draw,circle] (f) {$$}
;}\newcommand{\nodeg}{\node[draw,circle] (g) {$1$}
;}\newcommand{\nodeh}{\node[draw,circle] (h) {$$}
;}\newcommand{\nodei}{\node[draw,circle] (i) {$$}
;}\newcommand{\nodej}{\node[draw,circle] (j) {$$}
;}\newcommand{\nodeba}{\node[draw,circle] (ba) {$3$}
;}\newcommand{\nodebb}{\node[draw,circle] (bb) {$$}
;}\newcommand{\nodebc}{\node[draw,circle] (bc) {$$}
;}\newcommand{\nodebd}{\node[draw,circle,fill=black] (bd) {$$}
;}\newcommand{\nodebe}{\node[draw,circle] (be) {$2$}
;}\newcommand{\nodebf}{\node[draw,circle] (bf) {$1$}
;}\newcommand{\nodebg}{\node[draw,circle] (bg) {$$}
;}\newcommand{\nodebh}{\node[draw,circle] (bh) {$$}
;}\newcommand{\nodebi}{\node[draw,circle] (bi) {$$}
;}\newcommand{\nodebj}{\node[draw,circle] (bj) {$3$}
;}\newcommand{\nodeca}{\node[draw,circle] (ca) {$$}
;}\newcommand{\nodecb}{\node[draw,circle] (cb) {$$}
;}\newcommand{\nodecc}{\node[draw,circle] (cc) {$$}
;}\newcommand{\nodecd}{\node[draw,circle] (cd) {$3$}
;}\newcommand{\nodece}{\node[draw,circle] (ce) {$$}
;}\newcommand{\nodecf}{\node[draw,circle] (cf) {$$}
;}\newcommand{\nodecg}{\node[draw,circle,fill=blue] (cg) {$$}
;}\begin{tikzpicture}[auto]
\matrix[column sep=.3cm, row sep=.3cm,ampersand replacement=\&]{
         \&         \&         \&         \&         \&         \&         \&         \&         \& \nodea  \&         \&         \&         \&         \&         \&         \&         \&         \&         \\ 
         \&         \&         \&         \& \nodeb  \&         \&         \&         \&         \&         \&         \&         \&         \&         \& \nodebe \&         \&         \&         \&         \\ 
         \& \nodec  \&         \&         \& \nodeg  \&         \&         \& \nodeba \&         \&         \&         \& \nodebf \&         \&         \& \nodebj \&         \&         \& \nodecd \&         \\ 
 \noded  \& \nodee  \& \nodef  \& \nodeh  \& \nodei  \& \nodej  \& \nodebb \& \nodebc \& \nodebd \&         \& \nodebg \& \nodebh \& \nodebi \& \nodeca \& \nodecb \& \nodecc \& \nodece \& \nodecf \& \nodecg \\
};

\path[ultra thick, red] 
        (c) edge node [left]  {$2$} (d) edge node [right]  {\hspace{-.1cm}$1$} (e) edge node [right]  {$3$} (f)
	(g) edge node [left]  {$1$} (h) edge node [right]  {\hspace{-.1cm}$2$} (i) edge node [right]  {$3$} (j)
	(ba) edge node [left]  {$3$} (bb) edge node [right]  {\hspace{-.1cm}$1$} (bc) edge node [right]  {$2$} (bd)
	(b) edge node [near start, left]  {$1\quad$} (c) edge node [near start, right]  {\hspace{-.05cm}$3$} (g) edge node [near start, right]  {$\quad 2$} (ba)
	(bf) edge node [left]  {$1$} (bg) edge node [right]  {\hspace{-.1cm}$2$} (bh) edge node [right]  {$3$} (bi)
	(bj) edge node [left]  {$3$} (ca) edge node [right]  {\hspace{-.1cm}$2$} (cb) edge node [right]  {$1$} (cc)
	(cd) edge node [left]  {$1$} (ce) edge node [right]  {\hspace{-.1cm}$2$} (cf) edge node [right]  {$3$} (cg)
	(be) edge node [near start, left]  {$1\quad$} (bf) edge node [near start, right]  {\hspace{-.05cm}$2$} (bj) edge node [near start, right]  {$\quad 3$} (cd)
	(a) edge node [near start, left]  {$1\qquad$} (b) edge node [near start, right]  {$\qquad 2$} (be);
\end{tikzpicture}}
\end{center}
\caption{Graphical presentation of an elliptic map with $\RZ$ component action using the same labeling as in
Figure~\ref{figure.X}. The black leaf has coordinates $(1,2,2)$. Since it passes the constant maps 2,3,3 on its
way, it gets mapped to the leaf with coordinates $(2,3,3)$, denoted by the blue leaf.
\label{figure.RZ}}
\end{figure}

In light of the Prime Decomposition Theorem, there are three main cases for the component actions in $S_k$
of the elliptic maps on $T(n_0,\ldots,n_N)$. All of the next three statements have the following form. First note 
that composition of elliptic maps on a fixed tree with component action in a fixed pseudovariety is closed
under composition. Suppose that the component action $S_k$ is selected to be in the pseudovariety $\mathbf{V}$.
Then the pseudovariety generated by elliptic maps with component action in $\mathbf{V}$
(in this case divisors of elliptic maps) is determined and is denoted $\mathbf{PV}(\text{component in }\mathbf{V})$.
It is the semigroups of $\mathbf{PV}(\text{components in }\mathbf{V})$ on which we analyze their random walks:
\begin{enumerate}
\item
$S_k$ is in the \defn{pseudovariety $\RZ$} with $\mathbf{PV}(\text{component in }\RZ)$ which is delay semigroups
(see Section~\ref{subsection.delay}). In this case the component action consists only of constant maps.
If we label the branches from a vertex at level $k$ by $X_k=\{1,2,\ldots,n_k\}$, then we can also label the vertices 
at level $k$ by elements in $X_k$. The label $a \in X_k$ means the constant map that maps everything to $a$.
An example is given in Figure~\ref{figure.RZ}.
\item
$S_k$ is in the \defn{pseudovariety $\RZ^1$} with $\mathbf{PV}(\text{component in }\RZ^1)$
which is aperiodic semigroups (which means semigroups with  trivial subgroups).
In this case the component action consists of constant maps and 
the identity; the component monoids are aperiodic. If again the branches at level $k$ are labelled by
$X_k=\{1,2,\ldots,n_k\}$, then we can label the vertices by elements in $X_k \cup \{I\}$, where as before 
$a\in X_k$ denotes the constant map to $a$ and $I$ is the identity.
\item
$S_k$ is \defn{any finite group plus constant maps} and $\mathbf{PV}(\text{component in any finite 
group plus constant maps})$ is all finite semigroups.
In this case the vertices at level $k$ are labelled by elements in a finite
group $G$ which acts on the right on $X_k$ and elements in $X_k$ which give the constant maps. This yields 
a component semigroup with group of units in $G$ and kernel in $\RZ$.
\end{enumerate}

In this paper we will restrict to elliptic maps or wreath products with component actions in $\RZ$, that is 
constant maps  (without identity) to answer the question about resets. Future papers will deal with cases
2 and 3. 

\subsection{Delay pseudovariety}
\label{subsection.delay}

Let ${\bf D}$ be the pseudovariety of semigroups whose idempotents are right
zeroes, also called the 
\defn{delay pseudovariety}. The pseudovariety ${\bf D}$ can be
characterized (see~\cite[pg. 248]{RS:2009}) by 
\[
	{\bf D} = \bigcup_{k\ge 1} {\bf D}_k,
\]
where
\begin{equation}
	{\bf D}_k = [[ x_0 x_1 \cdots x_k = x_1 \cdots x_k]]\;,
\end{equation}
meaning that any $k+1$ elements $x_0,\ldots, x_k$ in a semigroup $S
\in {\bf D}_k$ satisfy the identity $x_0 x_1 \cdots x_k = x_1 \cdots x_k$.

The delay pseudovariety is also equal to $\RZnil$ defined as
\[
	\RZnil = \{ S \mid \text{$S/\ker(S)$ is nilpotent and
          $\ker(S)\in \RZ$}\}\;, 
\]
where we recall that $\RZ=[[xy=y]]$. A semigroup $N$ with zero is
nilpotent if $N^k= \{ 0 \} $ for some $k$, or 
in other words, $x_1 \cdots x_k=0$ in $N$. Thus, $S\in {\bf D}$ if and
only if $S$ satisfies the pseudoidentity  
$xy^{\omega} = y^{\omega}$, where $y^\omega$ is the unique idempotent
in $\langle y \rangle \le S$, or  
more succinctly
\[
	{\bf D} = [[xy^\omega = y^\omega]] = \RZnil \;.
\]
The pseudovariety ${\bf D}$ is also closed under semidirect products.
For all details see~\cite{RS:2009}.

\medskip

A semigroup $S$ is a \defn{subdirect product} of $S_1$ and $S_2$,
denoted $S \ll S_1 \times S_2$, if $S$ 
is a subsemigroup of $S_1 \times S_2$ mapping onto both $S_1$ and $S_2$ via the 
projections~\cite[pg. 34]{RS:2009}. More concretely, $S \ll S_1 \times S_2$ if and only if there
exist surmorphisms $\varphi_i\colon S \to S_i$ for $i = 1,2$, 
so that $\varphi_1$ and $\varphi_2$ separate points, that is, $s,t \in S$ with $s\neq t$ implies that 
$(s)\varphi_j \neq (t)\varphi_j$ for some $j \in \{1,2\}$.
The \defn{right letter mapping congruence} on a semigroup $S \in {\bf D}$
is defined by $s \sim t$ if $zs=zt$
for all $z\in \ker(S)$, that is, we identify two elements of $S$ if
they act the same on the right of $\ker(S)$. Therefore $\sim$ is the
kernel of the right Sch\"utzenberger representation of $S$ on
$\ker(S)$. We denote by 
$\RLM \colon S \twoheadrightarrow S$ the canonical morphism $s \mapsto
s/\sim$, and denote its image by $\RLM(S)$. (This definition agrees
with the definition given in~\cite[Section 4.6.2]{RS:2009}). 

From this it now follows that if $S \in {\bf D} = \RZnil$, then
\[
	S \ll S/\ker(S) \times \RLM(S).
\]
This can be observed by letting $\varphi_1 \colon S \to S/\ker(S)$ be the Rees quotient map, which maps 
$s\mapsto s$ if $s\not \in \ker(S)$ and collapses $\ker(S)$ to a single element.
Let $\varphi_2 \colon S \to \RLM(S)$ be the map $s\mapsto s/\sim$.
Hence $\varphi_2$ is injective on $\ker(S)$, so that $\varphi_1$ and $\varphi_2$ separate points.
In our applications, we only care about $\RLM(S)$. Note that a
semigroup $S\in {\bf D}$ is nilpotent if and only if
$\RLM(S)$ is the trivial semigroup $(0)$.

Observe that for $S,T\in {\bf D}$ we have $\ker(S),\ker(T)\in \RZ$ and
\begin{equation}
\begin{split}
 & \text{if} \quad S \twoheadrightarrow T \quad \text{then} \quad \RLM(S) \twoheadrightarrow \RLM(T)\\
 & \text{if} \quad S \twoheadrightarrow T \quad \text{then} \quad \ker(S) \twoheadrightarrow \ker(T)\\
 & \RLM(\RLM(S)) \cong \RLM(S).
\end{split}
\end{equation}
The proofs are not difficult and all details can be found in~\cite[Section 4.6.2]{RS:2009}.

\begin{Definition}
\label{definition.right congruence}
An equivalence relation $\tau$ on $\ker(S)$ is called a \defn{right congruence} if 
it preserves the right action of $S$ on $\ker(S)$,
that is, if $z \tau z'$ implies $(zs)\tau(z's)$ for all $z,z' \in
\ker(S)$ and $s \in S$. 
We denote by $\rc(\ker(S),S)$ (or by
$\rc(\ker(S))$ if $S$ is implicit) the set of all right
congruences on $\ker(S)$.
\end{Definition}

We consider $\rc(\ker(S))$ (partially) ordered by inclusion.
Since the intersection of right congruences on $\ker(S)$
is still a right congruence, $(\rc(\ker(S)),\subseteq)$ is a (complete)
$\wedge$-semilattice. Thus $(\rc(\ker(S)),\subseteq)$ is indeed a
(complete) lattice with the 
determined join, described by
$$\vee \Lambda = \bigcap\{ \rho \in \rc(\ker(S)) \mid \lambda
\subseteq \rho\mbox{ 
  for every }\lambda \in \Lambda\}$$
 for every $\Lambda \subseteq \rc(\ker(S))$.

It is routine to check each $\tau \in \rc(\ker(S),S)$ determines a
congruence $\oo{\tau}$ on $(\ker(S),\RLM(S))$ defined by 
$$(s\sim)\oo{\tau}(t\sim) \mbox{ if $(zs)\tau(zt)$ for every }z \in \ker(S),$$
where $s\sim$ denotes the equivalence class of $s \in S$ under the right letter mapping congruence $\sim$.
Since $S \in {\bf D}$, we have $\ker(S) \in \RZ$, and it follows easily that 
\beq
\label{ootau}
z\tau z'\mbox{ if and only if $(z\sim)\oo{\tau}(z'\sim)$ holds
for all }z,z' \in \ker(S).
\eeq
Thus right congruences on $\ker(S)$ and right letter mapping images of $S$
are the ``same thing''.

\subsection{Right zero component action}
\label{rzca}

In this section, we specialize the elliptic maps on rooted uniformly branching trees of Section~\ref{subsection.elliptic} 
to the constant component action. That is, we restrict ourselves to the case that the component action
$S_{\ell} \in \RZ=[[xy=y]]$ for all $0\le \ell \le N$.

\medskip

Let $F(g,k)$ be the semigroup generated by $A_g := \{a_1,a_2,\ldots,
a_g\}$ modulo all relations of the form
$$a_{i_0}a_{i_1}\ldots a_{i_k} = a_{i_1}\ldots a_{i_k}$$
for $i_0,\ldots,i_k \in \{ 1,\ldots,g\}$. This semigroup admits a
convenient normal form: we can identify $F(g,k)$
with $A^{\le k} \setminus \{ \varepsilon\}$, the set of all nonempty
words on $A$ of length at 
most $k$ (we denote the empty word by $\varepsilon$). Note that we may
define length of an element of $F(g,k)$ as the length of the
respective normal form in $A^{\le k} \setminus \{ \varepsilon \}$.

Given $u \in A^+$, let
$u\xi_k$ denote the suffix of length $k$ of $u$ if $|u| \geq k$ and
$u$ otherwise.
We define a binary operation $\circ$ on $A^{\leq k} \setminus \{
\varepsilon \}$ by $$u \circ v = (uv)\xi_k.$$
This binary operation on the normal forms corresponds to the product
of $F(g,k)$. For example in $F(2,3)$ with $A_2=\{a,b\}$ we have
$aba\cdot a = baa$, $aba \cdot bbb = bbb$, $b\cdot a = ba$ and so on.

It is immediate that
$F(g,k)$ satisfies the identity
\begin{equation}
\label{equation.N}
	x_0 x_1 \cdots x_k = x_1 \cdots x_k.
\end{equation}
Indeed, $F(g,k)$ is the \defn{free pro}-${\bf
  D}_k$ semigroup over $A$ (see \cite[Subsection 3.2.2]{RS:2009} for details
on free pro-$\V$ 
semigroups, for a pseudovariety $\V$). Since $F(g,k)$ is finite, it follows that
$F(g,k) \in {\bf D}$. Note that we 
can identify $\ker(F(g,k))$ with $A^k$, the set of all words on $A$ of
length $k$.   

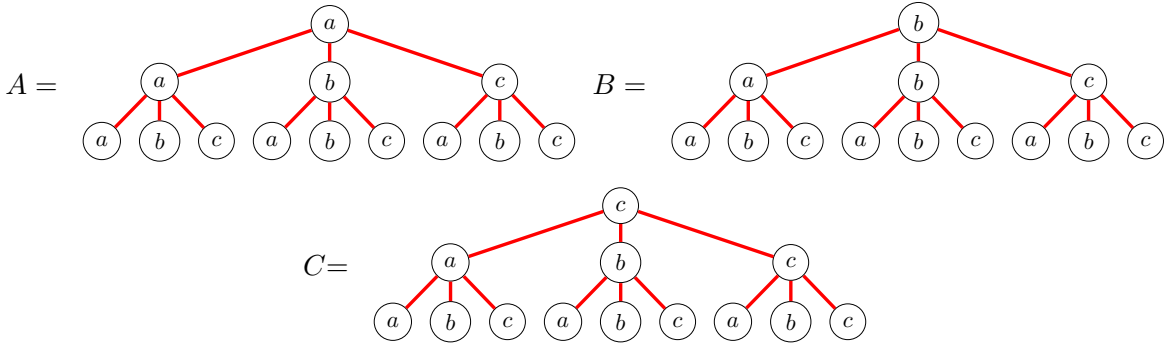
\begin{figure}[t]
\begin{center}
$A=$
\raisebox{-1cm}{\scalebox{0.8}
{ \newcommand{\nodea}{\node[draw,circle] (a) {$a$}
;}\newcommand{\nodeb}{\node[draw,circle] (b) {$a$}
;}\newcommand{\nodec}{\node[draw,circle] (c) {$a$}
;}\newcommand{\noded}{\node[draw,circle] (d) {$b$}
;}\newcommand{\nodee}{\node[draw,circle] (e) {$c$}
;}\newcommand{\nodef}{\node[draw,circle] (f) {$b$}
;}\newcommand{\nodeg}{\node[draw,circle] (g) {$a$}
;}\newcommand{\nodeh}{\node[draw,circle] (h) {$b$}
;}\newcommand{\nodei}{\node[draw,circle] (i) {$c$}
;}\newcommand{\nodej}{\node[draw,circle] (j) {$c$}
;}\newcommand{\nodeba}{\node[draw,circle] (ba) {$a$}
;}\newcommand{\nodebb}{\node[draw,circle] (bb) {$b$}
;}\newcommand{\nodebc}{\node[draw,circle] (bc) {$c$}
;}\begin{tikzpicture}[auto]
\matrix[column sep=.3cm, row sep=.3cm,ampersand replacement=\&]{
         \&         \&         \&         \& \nodea  \&         \&         \&         \&         \\ 
         \& \nodeb  \&         \&         \& \nodef  \&         \&         \& \nodej  \&         \\ 
 \nodec  \& \noded  \& \nodee  \& \nodeg  \& \nodeh  \& \nodei  \& \nodeba \& \nodebb \& \nodebc \\
};

\path[ultra thick, red] (b) edge (c) edge (d) edge (e)
	(f) edge (g) edge (h) edge (i)
	(j) edge (ba) edge (bb) edge (bc)
	(a) edge (b) edge (f) edge (j);
\end{tikzpicture}}}
$B=$
\raisebox{-1cm}{\scalebox{0.8}
{ \newcommand{\nodea}{\node[draw,circle] (a) {$b$}
;}\newcommand{\nodeb}{\node[draw,circle] (b) {$a$}
;}\newcommand{\nodec}{\node[draw,circle] (c) {$a$}
;}\newcommand{\noded}{\node[draw,circle] (d) {$b$}
;}\newcommand{\nodee}{\node[draw,circle] (e) {$c$}
;}\newcommand{\nodef}{\node[draw,circle] (f) {$b$}
;}\newcommand{\nodeg}{\node[draw,circle] (g) {$a$}
;}\newcommand{\nodeh}{\node[draw,circle] (h) {$b$}
;}\newcommand{\nodei}{\node[draw,circle] (i) {$c$}
;}\newcommand{\nodej}{\node[draw,circle] (j) {$c$}
;}\newcommand{\nodeba}{\node[draw,circle] (ba) {$a$}
;}\newcommand{\nodebb}{\node[draw,circle] (bb) {$b$}
;}\newcommand{\nodebc}{\node[draw,circle] (bc) {$c$}
;}\begin{tikzpicture}[auto]
\matrix[column sep=.3cm, row sep=.3cm,ampersand replacement=\&]{
         \&         \&         \&         \& \nodea  \&         \&         \&         \&         \\ 
         \& \nodeb  \&         \&         \& \nodef  \&         \&         \& \nodej  \&         \\ 
 \nodec  \& \noded  \& \nodee  \& \nodeg  \& \nodeh  \& \nodei  \& \nodeba \& \nodebb \& \nodebc \\
};

\path[ultra thick, red] (b) edge (c) edge (d) edge (e)
	(f) edge (g) edge (h) edge (i)
	(j) edge (ba) edge (bb) edge (bc)
	(a) edge (b) edge (f) edge (j);
\end{tikzpicture}}}
\end{center}

\begin{center}
$C$=
\raisebox{-1cm}{\scalebox{0.8}
{ \newcommand{\nodea}{\node[draw,circle] (a) {$c$}
;}\newcommand{\nodeb}{\node[draw,circle] (b) {$a$}
;}\newcommand{\nodec}{\node[draw,circle] (c) {$a$}
;}\newcommand{\noded}{\node[draw,circle] (d) {$b$}
;}\newcommand{\nodee}{\node[draw,circle] (e) {$c$}
;}\newcommand{\nodef}{\node[draw,circle] (f) {$b$}
;}\newcommand{\nodeg}{\node[draw,circle] (g) {$a$}
;}\newcommand{\nodeh}{\node[draw,circle] (h) {$b$}
;}\newcommand{\nodei}{\node[draw,circle] (i) {$c$}
;}\newcommand{\nodej}{\node[draw,circle] (j) {$c$}
;}\newcommand{\nodeba}{\node[draw,circle] (ba) {$a$}
;}\newcommand{\nodebb}{\node[draw,circle] (bb) {$b$}
;}\newcommand{\nodebc}{\node[draw,circle] (bc) {$c$}
;}\begin{tikzpicture}[auto]
\matrix[column sep=.3cm, row sep=.3cm,ampersand replacement=\&]{
         \&         \&         \&         \& \nodea  \&         \&         \&         \&         \\ 
         \& \nodeb  \&         \&         \& \nodef  \&         \&         \& \nodej  \&         \\ 
 \nodec  \& \noded  \& \nodee  \& \nodeg  \& \nodeh  \& \nodei  \& \nodeba \& \nodebb \& \nodebc \\
};

\path[ultra thick, red] (b) edge (c) edge (d) edge (e)
	(f) edge (g) edge (h) edge (i)
	(j) edge (ba) edge (bb) edge (bc)
	(a) edge (b) edge (f) edge (j);
\end{tikzpicture}}}
\end{center}
\caption{Generators for $F(3,3)$ on $T(3,3,3)$ with $A_3=\{a,b,c\}$.
\label{figure.tree}}
\end{figure}

It can also be interpreted in terms of elliptic maps on $T:=T(\underbrace{g,\ldots,g}_k)$ as follows. As in 
Section~\ref{subsection.elliptic}, we represent elliptic maps directly on the tree by denoting the component 
action on the vertices. Define the generators
$\varphi_1,\ldots,\varphi_g$ through trees of depth $k$ with $g$ branches at each level, 
where in level $1\le \ell \le k$ the vertices are labeled $a_1,\ldots,a_g$ from left to right. The $i$-th generator 
has label $a_i$ at level $0$. Since the vertices at level $k$ are not labeled, we will omit them for space reasons.
An example of the generators for $F(3,3)$ is given in Figure~\ref{figure.tree}.

\begin{figure}[t]
\begin{center}
$A\cdot B =$
\raisebox{-1cm}{\scalebox{0.8}
{ \newcommand{\nodea}{\node[draw,circle] (a) {$b$}
;}\newcommand{\nodeb}{\node[draw,circle] (b) {$a$}
;}\newcommand{\nodec}{\node[draw,circle] (c) {$a$}
;}\newcommand{\noded}{\node[draw,circle] (d) {$a$}
;}\newcommand{\nodee}{\node[draw,circle] (e) {$a$}
;}\newcommand{\nodef}{\node[draw,circle] (f) {$a$}
;}\newcommand{\nodeg}{\node[draw,circle] (g) {$b$}
;}\newcommand{\nodeh}{\node[draw,circle] (h) {$b$}
;}\newcommand{\nodei}{\node[draw,circle] (i) {$b$}
;}\newcommand{\nodej}{\node[draw,circle] (j) {$a$}
;}\newcommand{\nodeba}{\node[draw,circle] (ba) {$c$}
;}\newcommand{\nodebb}{\node[draw,circle] (bb) {$c$}
;}\newcommand{\nodebc}{\node[draw,circle] (bc) {$c$}
;}\begin{tikzpicture}[auto]
\matrix[column sep=.3cm, row sep=.3cm,ampersand replacement=\&]{
         \&         \&         \&         \& \nodea  \&         \&         \&         \&         \\ 
         \& \nodeb  \&         \&         \& \nodef  \&         \&         \& \nodej  \&         \\ 
 \nodec  \& \noded  \& \nodee  \& \nodeg  \& \nodeh  \& \nodei  \& \nodeba \& \nodebb \& \nodebc \\
};

\path[ultra thick, red] (b) edge (c) edge (d) edge (e)
	(f) edge (g) edge (h) edge (i)
	(j) edge (ba) edge (bb) edge (bc)
	(a) edge (b) edge (f) edge (j);
\end{tikzpicture}}}
\end{center}
\caption{Multiplication of elements $A$ and $B$ in $F(3,3)$. Note that the first two levels are constant precisely
as specified by $A$ and $B$.
\label{figure.mult}}
\end{figure}
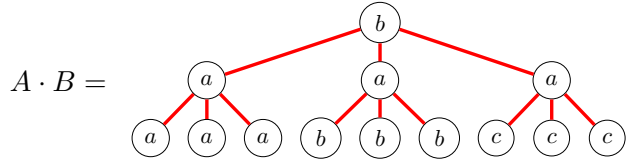

A label $a_i$ in a given vertex denotes the constant map to $a_i$. If we label the edges under
each vertex also $a_1,\ldots,a_g$ from left to right, then we can multiply generators on the labeled
tree as in Section~\ref{subsection.elliptic}. See Figure~\ref{figure.mult} for the product of $A$ and $B$ 
of Figure~\ref{figure.tree}. Using the notation $v_{j_1\ldots
  j_{k}}$ to denote the nodes below the root as in Subsection
\ref{subsection.elliptic}, we have $v_{j_1\ldots j_k}\p_i =
v_{ij_1\ldots j_{k-1}}$ and so 
$$v_{j_1\ldots j_k}\p_{i_{\ell-1}}\ldots \p_{i_{0}} = v_{i_0\ldots
  i_{\ell-1}j_1\ldots j_{k-\ell}}$$
for every $\ell \leq k$. In terms of component actions, this translates
into a tree with
$a_{i_0}$ on level 0, $a_{i_1}$ on all $g$ vertices of level 1, and in
general $a_{i_j}$ on all vertices of level $j$ for $0\le j<\ell$. It
follows easily from 
$$v_{j_1\ldots j_k}\p_{i_{k-1}}\ldots \p_{i_{0}} = v_{i_0\ldots
  i_{k-1}} = v_{j_1\ldots j_k}\p_{i_{k}}\ldots \p_{i_{0}}$$
that $\p_1,\ldots,\p_g$ generate a semigroup isomorphic to $F(g,k)$.

This gives a simple proof of Stiffler's Theorem~\cite{Stiffler:1973} (see also~\cite[Theorem 4.5.7, pg. 248]{RS:2009}).

\begin{T}[Stiffler]
The smallest pseudovariety containing the 2-element right zero semigroup that is closed under semidirect 
product (equivalently wreath or elliptic products) is ${\bf D}$.
\end{T}

\proof
As discussed in Section~\ref{subsection.delay}, ${\bf D}$ is a pseudovariety that is closed under 
semidirect product. 
By the arguments above, the free objects $F(g, k)$ are elliptic products with component action in $\RZ$ and since 
every member of ${\bf D}$ is a suromorphic image of an appropriate free one, the theorem is proved. \qed

In the sequel, we will be interested in the classification of right congruences on $\ker(F(g,k))\in \RZ$.

\section{$k$-reset graphs}
\label{regr}

$k$-reset graphs are finite state automata~\cite{Rhodes:2010} with the additional property that
strings of length $k$ are resets or constant maps. The formalism is such that the definitions in the profinite
case, when $k$ tends to infinity, is very similar. Let us now discuss the details.

Let $A$ be a finite nonempty alphabet.
An \defn{$A$-graph} is a structure of the
form $\Gamma = (Q,E)$, where:
\bi
\item
$Q$ is a finite nonempty set (vertex set);
\item
$E \subseteq Q \times A \times Q$ (edge set).
\ei

A \defn{nontrivial path} in an $A$-graph $\Gamma = (Q,E)$ is a finite
sequence of the form
$$q_0 \mapright{a_1} q_1 \mapright{a_2} \cdots \mapright{a_n} q_n$$
such that $(q_{i-1},a_i,q_i) \in E$ for $i = 1,\ldots,n$. Its label is
the word $a_1a_2 \cdots a_n \in A^+ = A^* \setminus \{\varepsilon\}$,
where $A^*$ is the set of words in the alphabet $A$ and $\varepsilon$ is
the empty word. A \defn{trivial path} is a formal expression of the form
$$q \mapright{\varepsilon} q.$$

An $A$-graph $\Gamma = (Q,E)$ is:
\bi
\item
\defn{deterministic} if
$$(p,a,q),(p,a,q') \in E \Rw q = q'$$
holds for all $p,q,q' \in Q$ and $a \in A$;
\item
\defn{complete} if 
$$\forall p \in Q\; \forall a \in A\; \exists q \in Q: (p,a,q) \in E;$$
\item
\defn{strongly connected} if, for all $p,q \in Q$, there exists a path
$p \mapright{u} q$ in $\Gamma$ for some $u \in A^*$.
\ei

If $\Gamma = (Q,E)$ is deterministic and complete, then $E$ induces a
function 
$$\begin{array}{rcl}
Q \times A&\to&Q\\
(q,a)&\mapsto&qa
\end{array}$$
defined by $(q,a,qa) \in E$. Conversely, every such function defines a
deterministic complete $A$-graph. Moreover, we
can extend the function $Q \times A \to Q$ to a function $Q \times
A^* \to Q$ as follows: given $q \in Q$ and $u \in A^*$, $qu$ is the
unique vertex such that there exists a path
$$q \mapright{u} qu$$ 
in $\Gamma$. This function is called the \defn{transition function} of
$\Gamma$. 

Let $\Gamma = (Q,E)$ and $\Gamma' = (Q',E')$ be $A$-graphs. A \defn{morphism} 
$\p:\Gamma \to \Gamma'$ is a function $\p:Q \to Q'$ such that 
$$(p,a,q) \in E \Rw (p\p,a,q\p) \in E'.$$
If $\p$ is bijective and $\p\inv$ is also a morphism, we say that $\p$
is an \defn{isomorphism}. In this case we write $\Gamma \cong \Gamma'$.

Given $A$-graphs $\Gamma,\Gamma'$, we write $\Gamma \leq \Gamma'$ if
there exists a morphism $\Gamma \to \Gamma'$. This is clearly a
reflexive and transitive relation, hence a preorder on the class of
all $A$-graphs. Technically, this is not a partial order, but we have
the following remark:

\bl
\label{paor}
Let $A$ be a finite nonempty alphabet and let $\Gamma,\Gamma'$ be
strongly connected deterministic complete $A$-graphs such that $\Gamma
\leq \Gamma' \leq \Gamma$. Then $\Gamma 
\cong \Gamma'$.   
\el

\proof
Let $\p:\Gamma \to \Gamma'$ and $\p':\Gamma' \to \Gamma$ be morphisms. 
Write $\Gamma = (Q,E)$ and $\Gamma' = (Q',E')$. Fix some $q_0 \in Q$
and take $q' \in Q'$. Since $\Gamma'$ is strongly connected, there
exists some path $q_0\p \mapright{u} q'$ in $\Gamma'$ for some $u \in
A^*$. Since $\Gamma$ is complete, there
exists some path $q_0 \mapright{u} q$ in $\Gamma$ for some $q \in
Q$. It follows from $\p$ being a morphism that there
exists a path $q_0\p \mapright{u} q\p$ in $\Gamma'$. Since $\Gamma'$
is deterministic, we get $q' = q\p$, hence $\p$ is onto and so $|Q'|
\leq |Q|$. By symmetry, we get $|Q'| = |Q|$, thus $\p$ is bijective.

It remains to be proved that $\p\inv$ is a morphism. Assume that
$(p\p,a,q\p) \in E'$ for some $p,q \in Q$ and $a \in A$. Since
$\Gamma$ is complete, there exists some $(p,a,r) \in E$. Since $\p$ is
a morphism, we get $(p\p,a,r\p) \in E'$. Now $\Gamma'$ being
deterministic yields $q\p = r\p$, and so $q = r$ since $\p$ is
bijective. Therefore $(p,a,q) \in E$ and so $\p\inv$ is a morphism as
required.
\qed

We say that $u \in A^*$ is a \defn{reset word} for the deterministic
and complete $A$-graph $\Gamma = (Q,E)$ if
$|Qu| = 1$. This is equivalent to say that all paths labeled by $u$
end at the same vertex. Let $\res(\Gamma)$ denote the set of all reset
words for $\Gamma$. For every $k \in \N$, let 
$$\res_k(\Gamma) = \res(\Gamma) \cap A^k.$$
We say that $\Gamma$ is a \defn{$k$-reset graph} if $\res_k(\Gamma) =
A^k$. We denote by $\rg_k(A)$ the class of all strongly connected
deterministic complete $k$-reset $A$-graphs.

Given $\Gamma \in \rg_k(A)$, let $[\Gamma]$ denote the isomorphism
class of $\Gamma$. Let
$$
	\rg_k(A) /\cong \; = \{ [\Gamma] \mid \Gamma \in \rg_k(A) \}.
$$
Given $\Gamma, \Gamma' \in \rg_k(A)$, write
$$
	[\Gamma] \leq [\Gamma'] \mbox{ if } \Gamma \leq \Gamma'.
$$
It is immediate that $\leq$ is a well-defined preorder on $\rg_k(A)
/\cong$. Moreover, it follows from Lemma~\ref{paor} that:

\bc
\label{porg}
Let $A$ be a finite nonempty alphabet and let $k \geq 1$. Then $\leq$
is a partial order on $\rg_k(A) /\cong$.
\ec

\section{Semaphore codes}
\label{section.semaphore codes}

A detailed discussion on semaphore codes can be found in~\cite[Chapter 3.4]{BPR:2010}.

Let $A$ be a finite alphabet. We define three partial orders on $A^*$ by
\bi
\item
$u \leq_p v$ if $v \in uA^*$,
\item
$u \leq_s v$ if $v \in A^*u$,
\item
$u \leq_f v$ if $v \in A^*uA^*$.
\ei
We refer to them as the \defn{prefix order}, the 
\defn{suffix order} and
the \defn{factor order} on $A^*$.

If $X \subset A^*$ is a nonempty antichain with respect to $\leq_p$
(respectively $\leq_s$, $\leq_f$), it is said to be a \defn{prefix code}
(respectively \defn{suffix code}, \defn{infix code}). 
Note that our notions differ slightly from the standard notions since
we admit $\{ \varepsilon \}$ to be a code of all three types!

Given an ideal $I \unlhd A^*$, let $I\beta$
denote the subset of elements of $I$ wich are minimal with respect to
$\leq_f$. Then $I = A^*(I\beta)A^*$ and $I\beta \subseteq B$ whenever $B
\subseteq A^*$ satisfies $I = A^*BA^*$. We say that $I\beta$ is the
\defn{basis} of $I$. Clearly, the correspondences
$$I \mapsto I\beta,\quad C \mapsto A^*CA^*$$
establish mutually inverse bijections between the set of all ideals of
$A^*$ and the set of all infix codes on $A$. 

We say that $L \subseteq A^*$ is a \defn{left ideal} if $L \neq \emptyset$ and
$A^*L \subseteq L$. We write then $L \unlhd_{\ell} A^*$.
Given $L \unlhd_{\ell} A^*$, let $L\beta_{\ell}$
denote the subset of elements of $L$ wich are minimal with respect to
$\leq_s$. Then $L = A^*(L\beta_{\ell})$ and $L\beta \subseteq B$ whenever $B
\subseteq A^*$ satisfies $L = A^*B$. We say that $L\beta_{\ell}$ is the
\defn{left basis} of $L$. 
Clearly, the correspondences
$$L \mapsto L\beta_{\ell},\quad S \mapsto A^*S$$
establish mutually inverse bijections between the set of all left ideals of
$A^*$ and the set of all suffix codes on $A$. 

Similarly, $R \subseteq A^*$ is a \defn{right ideal} if $R \neq \emptyset$ and
$RA^* \subseteq R$. We write then $R \unlhd_{r} A^*$.

We relate now ideals to semaphore codes. The definition we use is actually 
the left-right dual of the classical definition in~\cite[Section 3.5]{BPR:2010}, but we 
shall call them semaphores codes for simplification. We also admit $\emptyset$ 
and $\{ \varepsilon \}$ as (semaphore) codes, but this generalization is compatible with 
the relevant results from~\cite{BPR:2010}.

A \defn{semaphore code} on the alphabet $A$ is a language of the form
$$XA^* \setminus A^+XA^*,$$
for some $X \subseteq A^*$. If $X \neq \emptyset$, then $XA^*
\setminus A^+XA^*$ is a maximal suffix code (with respect to
inclusion) by~\cite[Proposition 3.5.1]{BPR:2010}. 
Now~\cite[Proposition 3.5.4]{BPR:2010} provides an alternative characterization of
semaphore codes:

\bl
\label{alch}
{\rm \cite[Proposition 3.5.4]{BPR:2010}}
For every $S \subseteq A^*$, the following conditions are equivalent:
\bi
\item[(i)]
$S$ is a semaphore code;
\item[(ii)]
$S$ is a suffix code and $SA \subseteq A^*S$. 
\ei
\el

Let $\sem(A)$ denote the set of all semaphore codes on the alphabet
$A$. We define a partial order $\leq$ on $\sem(A)$ by $S \leq S'$ if
$A^*S \leq A^*S'$. 

\begin{Example}
Let $A=\{a,b\}$ and $X=\{b\}$. Then the semaphore code is infinite
\[
	S=XA^* \setminus A^+ X A^*= \{b,ba,ba^2,ba^3,\ldots\} = ba^*.
\]
If on the other hand $A=\{a,b\}$ and $X=\{a^2,ab,b^2\}$, then the semaphore code is finite
\[
	S=XA^* \setminus A^+ X A^* = \{a^2,ab,b^2,aba,b^2a\}.
\]
\end{Example}

We denote by
$\I(A)$ (respectively $\L(A), \R(A)$) the set of all ideals (respectively
left ideals, right ideals) of $A^*$. 
If we  order $\I(A)$ (or $\L(A)$ or $\R(A)$) by inclusion,
we get a complete (distributive) lattice where meet and join are given
by intersection and union.
The top element is $A^*$ and the bottom element is $\emptyset$.
We can now prove the following.

\bp
\label{semaphore}
Let $A$ be a finite nonempty alphabet. Then
$$\begin{array}{rcl}
\Phi \colon (\I(A),\subseteq)&\to&({\rm Sem}(A),\leq)\\
I&\mapsto&I\beta_{\ell}
\end{array}
\quad\mbox{and}\quad
\begin{array}{rcl}
\Psi \colon ({\rm Sem}(A)\leq)&\to&(\I(A),\subseteq)\\ 
S&\mapsto&A^*S
\end{array}$$
are mutually inverse lattice isomorphisms.
\ep

\proof
Let $I \in \I(A)$. Then $I\beta_{\ell}$ is clearly a suffix code. Since
$(I\beta_{\ell})A \subseteq I = A^*(I\beta_{\ell})$, then
$I\beta_{\ell} \in \sem(A)$ by Lemma \ref{alch}
and $\Phi$ is well-defined.

On the other hand, given $S \in \sem(A)$, it is clear that $A^*S
\unlhd_{\ell} A^*$. Now $SA \subseteq A^*S$ by Lemma~\ref{alch}, hence $A^*S$ is
actually an ideal of $A^*$ and so $\Psi$ is also well-defined.

Now $I\Phi\Psi = A^*(I\beta_{\ell}) = I$ and $S\Psi\Phi =
(A^*S)\beta_{\ell} = S$ follows easily from $S$ being a suffix code,
hence $\Phi$ and $\Psi$ are mutually inverse bijections. Since $S \leq
S'$ if and only if $S\Psi \subseteq S'\Psi$ holds for all $S,S' \in
\sem(A)$, $\Phi$ and $\Psi$ are actually mutually inverse poset
isomorphisms. Since $(\I(A),\subseteq)$ is a lattice, so is
$(\sem(A),\leq)$ and so $\Phi$ and $\Psi$ are lattice isomorphisms.
\qed  

As we will see in Section~\ref{section.special right congruences}, semaphore codes
are related to special right congruences.

\section{Right congruences on the minimal ideal of $F(g,k)$}
\label{rcmi}

Now fix a nonempty alphabet $A = \{ a_1,\ldots,a_g \}$ and a positive
integer $k$. We remarked in Subsection~\ref{rzca} that $A^{\leq k}
\setminus \{ \varepsilon \}$ is a set of normal forms for $F(g,k)$,
the free pro-$\D_k$ semigroup on the set $A = \{ a_1,\ldots,a_g
\}$. Moreover, we can identify $A^k$ with $\ker(F(g,k))$.
Since $F(g,k)$ is generated by $A$, right congruences on $A^k$ can be
described as equivalence relations $\rho$ satisfying
$$u\rho v \Rw (u\circ a)\rho (v\circ a)$$
for every $a \in A$, or equivalently,
$$u\rho v \Rw ((ua)\xi_k)\rho ((va)\xi_k)$$
for every $a \in A$. 

Given $R \subseteq A^k \times A^k$, we denote by $R^{\sharp}$ the right
congruence on $A^k$ generated by $R$, i.e. the intersection of all right
congruences on $A^k$ containing $R$.
Let $u,v \in A^k$. Then $(u,v) \in
R^{\sharp}$ if and only if there exists some finite sequence $w_0,
\ldots, w_n \in A^k$ $(n \geq 0)$ such that:
\bi
\item
$w_0 = u$ and $w_n = v$;
\item
for every $i = 1,\ldots,n$, there exist $(r_i,s_i) \in R$ and $x_i \in
A^*$ such that $\{ w_{i-1},w_i \} = \{ r_i\circ x_i, s_i\circ x_i \}$.
\ei
It is easy to see that
$$\vee \Lambda = (\cup \Lambda)^{\sharp}$$
for every $\Lambda \subseteq \rc(A^k)$.

We now relate right congruences on $A^k$ with the $k$-reset graphs
introduced in Section~\ref{regr}. 

Given $\rho \in \rc(A^k)$, the \defn{Cayley graph} of $\rho$ is the
$A$-graph $\cay(\rho) = (A^k/\rho,E)$ defined by
$$E = \{ (u\rho,a,(u\circ a)\rho) \mid u \in A^k, \; a \in A \},$$
where $u\rho$ denotes the congruence class of $u$.
In particular, if $\rho$ is the identity relation, then $\cay(\rho)$
is a $k$-dimensional \defn{De Bruijn graph} on $|A|$ symbols.

Given $\Gamma = (Q,E) \in \rg_k(A)$, let $\zeta_{\Gamma}$ be the
equivalence relation on $A^k$ defined by
$$u \zeta_{\Gamma} v \mbox{ if } Qu = Qv.$$
Note that 
\beq
\label{isom3}
Q((ua)\xi_k) = Qua
\eeq
holds for all $u \in A^k$ and $a \in A$. Indeed, since $Qua \subseteq
Q((ua)\xi_k)$ and $(ua)\xi_k$ is a reset word, we must have equality
and (\ref{isom3}) holds.

\bp
\label{isom}
Let $A$ be a finite nonempty alphabet and $k \geq 1$. Then
$$\begin{array}{rcl}
\Phi \colon ({\rm RC}(A^k),\subseteq)&\to&({\rm RG}_k(A)/\cong,\leq)\\
\rho&\mapsto&[{\rm Cay}(\rho)]
\end{array}
\quad\mbox{and}\quad
\begin{array}{rcl}
\Psi \colon ({\rm RG}_k(A)/\cong,\leq)&\to&({\rm RC}(A^k),\subseteq)\\ 
{[\Gamma]}&\mapsto&\zeta_{\Gamma}
\end{array}$$
are mutually inverse lattice isomorphisms.
\ep

\proof
Let $\rho \in \rc(A^k)$. It follows from the definition that
$\cay(\rho)$ is deterministic and complete. For all $u,v \in A^k$, we
have $u \circ v = v$, hence there exists a path
$$u\rho \mapright{v} (u \circ v)\rho = v\rho$$
in $\cay(\rho)$. It follows that $\cay(\rho)$ is strongly connected
and $A^k \subseteq \res_k(\cay(\rho))$, thus $\cay(\rho) \in \rg_k(A)$
and $\Phi$ is well-defined.

On the other hand, it is clear that $[\Gamma]\Psi$ does not depend on
the chosen representative for the isomorphism class $[\Gamma]$.

Let $\Gamma \in \rg_k(A)$. Let $(u,v) \in \zeta_{\Gamma}$ and $a \in
A$. Then $Qu = Qv$ implies $Qua = Qva$ and therefore $(u
\circ a,v \circ a) \in \zeta_{\Gamma}$ in view of (\ref{isom3}). Thus
$\zeta_{\Gamma} \in 
\rc(A^k)$ and so $\Psi$ is well-defined.

Let $\rho \in \rc(A^k)$ and write $\rho' = \zeta_{\cay(\rho)}$. If
$Q = A^k/\rho$ is the vertex set of $\cay(\rho)$, then $Qu = \{ u\rho \}$
for every $u \in A^k$. Hence
$$u\rho'v \iff Qu = Qv \iff u\rho = v\rho$$
and so $\Phi\Psi = 1$.

Conversely, let $\Gamma = (Q,E) \in \rg_k(A)$ and let $\Gamma' =
\cay(\zeta_{\Gamma})$. We show that
\beq
\label{isom1}
\forall q \in Q \; \exists u_q \in A^k : Qu_q = \{ q \}.
\eeq
We may assume that $|Q| > 1$. Since $\Gamma$ is strongly connected, it
follows that there exists a loop $q \mapright{w} q$ in $\Gamma$ with
$w \neq \varepsilon$. Replacing $w$ by a proper power if necessary, we
may assume that $|w| \geq k$. Hence there exists some $u_q \in A^k$ such
that $q \in Qu_q$. Since $u_q$ is necessarily a reset word, we get $Qu_q =
\{ q \}$ and so (\ref{isom1}) holds.

We define a mapping
$$\begin{array}{rcl}
\theta \colon Q&\to&A^k/\zeta_{\Gamma}.\\
q&\mapsto&u_q\zeta_{\Gamma}
\end{array}$$
Note that 
\beq
\label{isom2}
Qu = Qv \iff u\zeta_{\Gamma} = v\zeta_{\Gamma}
\eeq
holds for all $u,v \in A^k$, hence $\theta$ is well-defined and
one-to-one. Since $\Gamma$ is a $k$-reset graph, $\theta$ is also
onto. We show that $\theta$ is an isomorphism from $\Gamma$ onto
$\cay(\zeta_{\Gamma})$. 

Assume that $(p,a,q) \in E$. By (\ref{isom3}), we get
$$Q(u_p \circ a) = Qu_pa = pa = q = Qu_q.$$
Hence $u_q\zeta_{\Gamma} = (u_p \circ a)\zeta_{\Gamma}$ and so
there exists an edge $u_p\zeta_{\Gamma} \mapright{a}
u_q\zeta_{\Gamma}$ in $\cay(\zeta_{\Gamma})$.

Conversely, assume that $u_p\zeta_{\Gamma} \mapright{a}
u_q\zeta_{\Gamma}$ is an edge of $\cay(\zeta_{\Gamma})$. Then 
$u_q\zeta_{\Gamma} = (u_p \circ a)\zeta_{\Gamma}$ and so 
$$q = Qu_q = Q(u_p \circ a) = Qu_pa = pa$$ 
by (\ref{isom2}) and (\ref{isom3}). Thus $(p,a,q) \in E$ and so
$\theta \colon \Gamma \to \cay(\zeta_{\Gamma})$ is an isomorphism. Therefore
$\Psi\Phi = 1$ and so $\Phi$ and $\Psi$ are mutually inverse
bijections.

Let $\rho,\rho' \in \rc(A^k)$ with $\rho \subseteq \rho'$. Then
$$\begin{array}{rcl}
\theta \colon A^k/\rho&\to&A^k/\rho'\\
u\rho&\mapsto&u\rho'
\end{array}$$
is a well-defined surjective map. If $u\rho \mapright{a} (u\circ
a)\rho$ is an edge of $\cay(\rho)$, then $u\rho' \mapright{a} (u\circ
a)\rho'$ is an edge of $\cay(\rho)'$, hence $\theta$ is a morphism
from $\cay(\rho)$ to $\cay(\rho')$ and so $\cay(\rho) \leq
\cay(\rho')$. Thus $[\cay(\rho)] \leq
[\cay(\rho')]$ and so $\Phi$ is order-preserving.

Let $\Gamma,\Gamma' \in \rg_k(A)$ be such that $[\Gamma] \leq
[\Gamma']$. Then there exists a morphism $\theta \colon \Gamma \to
\Gamma'$. Write $\Gamma = (Q,E)$ and $\Gamma' = (Q',E')$. Suppose that
$(u,v) \in \zeta_{\Gamma}$. Then $Qu = Qv = \{q\}$ for some $q \in
Q$. Hence $q\theta \in Q'u \cap Q'v$. Since $\Gamma'$ is a $k$-reset
graph, we get $Q'u = \{q\theta \} = Q'v$ and so $(u,v) \in
\zeta_{\Gamma'}$. Therefore $\Psi$ is order-preserving.

Since $\Phi$ and $\Psi$ are mutually inverse order-preserving
mappings, they are isomorphisms of posets. Since
$(\rc(A^k),\subseteq)$ is a lattice, then $(\rg_k(A),\leq)$ is also a
lattice, and $\Phi$ and $\Psi$ are mutually inverse lattice
isomorphisms.
\qed

\section{Lattice-theoretic properties}
\label{lrc}

We discuss in this section the lattice-theoretic properties of the
lattice $\rc(A^k)$. 

We recall some well-known notions from lattice theory. Let $L$ be a
(finite) lattice with bottom element $B$ and top element $T$. Given
$a,b \in L$, we say that $b$ \defn{covers} $a$ if $a < b$ and there is
no $c \in L$ such that $a < c < b$. If $a$ covers the bottom $B$, we
say that $a$ is an \defn{atom}.

The lattice $L$ is said to be:
\bi
\item
\defn{modular} if it has no sublattice of the
form
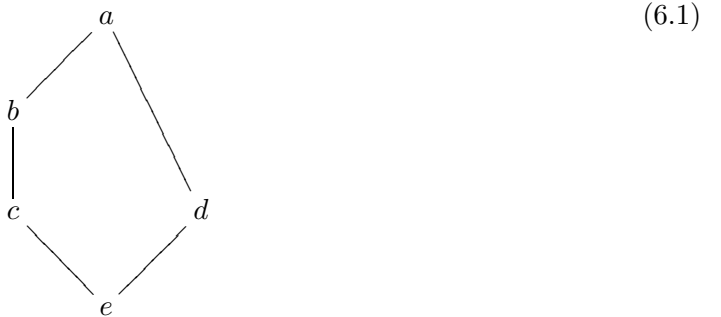
\beq
\label{modl}
\xymatrix{
& a \ar@{-}[dl] \ar@{-}[ddr] & \\
b \ar@{-}[d] && \\
c \ar@{-}[dr] && d \ar@{-}[dl] \\
& e &
}
\eeq
\item
\defn{semimodular} if it has no sublattice of the 
form (\ref{modl}) with $d$ covering $e$;
\item
\defn{atomistic} if every element of $L$ 
is a join of atoms ($B$ being the join of the empty set).
\ei

\bp
\label{lek}
Let $A$ be a nonempty set and $k \geq 1$. Then ${\rm RC}(A^k)$ is semimodular.
\ep

\proof
It suffices to show that $\rc(A^k)$ has no sublattice of the
form
$$
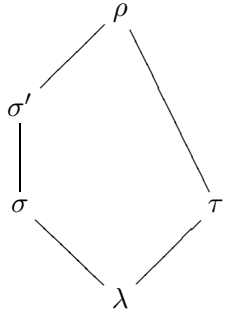
\xymatrix{
& \rho \ar@{-}[dl] \ar@{-}[ddr] & \\
\sigma' \ar@{-}[d] && \\
\sigma \ar@{-}[dr] && \tau \ar@{-}[dl] \\
& \lambda &
}$$
with $\tau$ covering $\lambda$ in $\rc(A^k)$.

Suppose it does.
Given $x,y \in A^*$, let $\lcs(x,y)$ denote the longest common suffix
of $x$ and $y$. If $x,y \in A^k$ are distinct, then $|\lcs(x,y)| <
k$ and so
\beq
\label{lek1}
|\lcs(x\circ a,y\circ a)| > |\lcs(x,y)|
\eeq
for every $a \in A$. 

Let $(u,v) \in \tau \setminus \lambda$ with $|\lcs(u,v)|$ maximal. For
every $a \in A$, we have  
$$(u\circ
a,v\circ a) \in \{ (u,v)\}^{\sharp} \subseteq \tau.$$ 
In view of (\ref{lek1}), and by maximality of $|\lcs(u,v)|$, we get 
\beq
\label{lek3}
(u\circ a,v\circ a) \in \lambda.
\eeq
Note also that
$$\lambda \subset (\lambda \cup \{ (u,v)\})^{\sharp} \subseteq \tau$$
yields 
\beq
\label{lek2}
\tau = (\lambda \cup \{ (u,v)\})^{\sharp}
\eeq
since $\tau$ covers $\lambda$. 

Let $(y,z) \in \sigma' \setminus \sigma$. Then (\ref{lek2}) yields
$$(y,z) \in \rho = (\sigma \vee \tau) = (\sigma \cup (\lambda \cup \{
(u,v)\})^{\sharp})^{\sharp} = (\sigma \cup \{
(u,v)\} )^{\sharp}$$
and so there exists some finite sequence $w_0,
\ldots, w_n \in A^k$ such that:
\bi
\item
$w_0 = y$ and $w_n = z$;
\item
for every $i = 1,\ldots,n$, there exist $(r_i,s_i) \in \sigma \cup \{
(u,v)\}$ and $x_i \in
A^*$ such that $\{ w_{i-1},w_i \} = \{ r_i\circ x_i, s_i\circ x_i \}$.
\ei
Now by (\ref{lek3}) we may assume that $x_i = \varepsilon$ whenever $(r_i,s_i) =
(u,v)$. Since we may assume that the $w_i$ are all distinct, the
relation $(u,v)$ is used at most once, indeed exactly once since
$(y,z) \notin \sigma$ and $(r_i,s_i) \in \sigma$ implies $(r_i\circ
x_i,s_i\circ x_i) \in \sigma$. We may assume without loss of
generality that $u = w_{j-1}$ and $v = w_j$ for some $j \in \{
1,\ldots, n\}$. Hence
$$y = w_0 \,\sigma \, w_{j-1} = u,\quad v = w_j\, \sigma\, w_n = z$$
and so 
$$u = w_{j-1}\, \sigma'\, y \, \sigma'\, z \, \sigma'\, w_j = v.$$
It follows that $\lambda \cup \{ (u,v)\} \subseteq \sigma'$. By
(\ref{lek2}), we get $\tau \subseteq \sigma'$, a contradiction.
Therefore $\rc(A^k)$ is semimodular.
\qed

Since a semimodular lattice of finite height (i.e. the length of chains is
bounded) satisfies the Jordan-Dedekind condition
(i.e. all maximal chains have the same length), we immediately obtain:

\bc
\label{lekjd}
Let $A$ be a nonempty set and $k \geq 1$. Then ${\rm RC}(A^k)$
satisfies the Jordan-Dedekind condition.
\ec

We show next that we cannot replace semimodular by modular in
Proposition~\ref{lek}. 

\bp
\label{lekmod}
Let $k \geq 1$ and let $A$ be a set with $|A| \geq 4$. Then ${\rm
  RC}(A^k)$ is  not modular.
\ep

\proof
Let $a,b,c,d \in A$ be distinct. Let $\lambda$ be the identity
relation on $A^k$ and let
\bi
\item[]
$\sigma = \lambda \cup \{ a^k, ba^{k-1} \}^2$;
\item[]
$\sigma' = \lambda \cup \{ a^k, ba^{k-1} \}^2 \cup \{ ca^{k-1}, da^{k-1} \}^2$;
\item[]
$\tau = \lambda \cup \{ a^k, da^{k-1} \}^2 \cup \{ ba^{k-1}, ca^{k-1} \}^2$;
\item[]
$\rho = \lambda \cup \{ a^k, ba^{k-1}, ca^{k-1}, da^{k-1} \}^2$.
\ei
It is routine to check that all the above relations are right
congruences on $A^k$. Moreover, 
$$\lambda \subset \sigma \subset \sigma' \subset \rho, \quad
\lambda \subset \tau \subset \rho,$$ 
$$\sigma' \cap \tau = \lambda, \quad (\sigma \vee \tau) = \rho,$$
hence  
$$
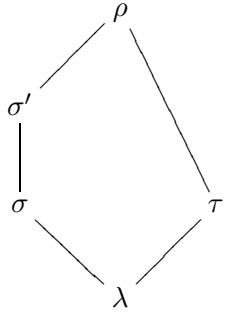
\xymatrix{
& \rho \ar@{-}[dl] \ar@{-}[ddr] & \\
\sigma' \ar@{-}[d] && \\
\sigma \ar@{-}[dr] && \tau \ar@{-}[dl] \\
& \lambda &
}$$
is a sublattice of $\rc(A^k)$ and so $\rc(A^k)$ is not
modular.
\qed

We can also show that $\rc(A^k)$ can only be atomistic in
trivial cases:

\bp
\label{lekat}
Let $k \geq 2$ and let $A$ be a set with $|A| \geq 2$. Then ${\rm
  RC}(A^k)$ is not atomistic.
\ep

\proof
Let $\lambda$ be the identity
relation on $A^k$. Let $a,b \in A$ be distinct and let 
\bi
\item[]
$\sigma = \lambda \cup \{ a^k, b^2a^{k-2}, ba^{k-1} \}^2
\cup \{ a^{k-1}b, ba^{k-2}b \}^2$;
\item[]
$\tau = \lambda \cup \{ a^k, ba^{k-1} \}^2 \cup \{ a^{k-1}b, ba^{k-2}b \}^2$.
\ei
It is routine to check that $\sigma, \tau \in
\rc(A^k)$. Moreover, $\lambda \subset \tau \subset \sigma$. 
We show that
\beq
\label{lekat1}
\sigma = \{ (xa^{k-1}, b^2a^{k-2}) \}^{\sharp}
\eeq
for every $x \in \{ a,b\}$. Indeed, let $\eta = \{ (xa^{k-1},
b^2a^{k-2}) \}^{\sharp}$. Then $(xa^{k-1}, b^2a^{k-2}) \in \eta$
yields $(a^{k}, ba^{k-1}) \in \eta$ and so $\{ a^k, b^2a^{k-2},
ba^{k-1} \}^2 \subseteq \eta$. Finally, $(xa^{k-1}, b^2a^{k-2}) \in
\eta$ yields $(a^{k-1}b, ba^{k-2}b) \in \eta$ and so
$$\sigma \subseteq \{ (xa^{k-1}, b^2a^{k-2}) \}^{\sharp}.$$
Since $(xa^{k-1}, b^2a^{k-2}) \in \sigma$ for $x \in \{ a,b\}$, (\ref{lekat1})
holds. 

Now we
claim that $\tau$ is the unique element of $\rc(A^k)$ covered
by $\sigma$. Indeed, assume that $\rho \subset \sigma$. In view of
(\ref{lekat1}), we have $(a^{k}, b^2a^{k-2}) \notin \rho$ and
$(ba^{k-1}, b^2a^{k-2}) \notin \rho$. Hence $\rho \subseteq
\tau$. Since $\sigma$ is not an atom, it follows that 
$$\alpha \leq \sigma \mbox{ if and only if } \alpha \leq \tau$$
for every atom $\alpha$ of $\rc(A^k)$. Thus $\sigma$ 
cannot be expressed as a join of atoms and so $\rc(A^k)$ is not
atomistic.
\qed

\section{Special right congruences on $A^k$}
\label{section.special right congruences}

To avoid trivial cases, we assume throughout this section that $A$ is
a finite alphabet containing at least two elements. We define
$$\I_k(A) = \{ I \unlhd A^* \mid A^k \subset I\},$$
$$\L_k(A) = \{ L \unlhd_{\ell} A^* \mid A^k \subset L\}.$$
If we order $\I_k(A)$ (or $\L_k(A)$) by inclusion, we get a finite
(distributive) lattice where meet and join are given by
$$(I \wedge J) = I \cap J,\quad (I \vee J) = I \cup J.$$
The top element is $A^*$ and the bottom element is $A^kA^*$.

Given $L \in \L_k(A)$, we define a relation $\tau_L$ on $A^k$ by:
$$u \tau_L v \; \mbox{ if $u$ and $v$ have a common suffix in }L.$$

\bl
\label{irco}
Let $L \in \L_k(A)$. Then $\tau_L$ is an equivalence relation on $A^k$.
\el

\proof
It is immediate that $\tau_L$ is symmetric. Since $A^k \subseteq
L$, it is reflexive. Assume now that $u,v,w \in A^k$ and $x,y \in L$
are such that $x \leq_s u,v$ and $y \leq_s v,w$. Since $x$ and $y$ are
both suffixes of $v$, one of them is a suffix of the other. Hence
either $x \leq_s u,w$ or $y \leq_s u,w$. Therefore $\tau_L$ is
transitive.
\qed

Being a right congruence turns out to be a special case:

\bp
\label{nsrc}
Let $L \in \L_k(A)$. Then the following conditions are equivalent:
\bi
\item[(i)] $\tau_L \in {\rm RC}(A^k)$;
\item[(ii)] $L \in \I_k(A)$;
\item[(iii)] $(L\beta_{\ell})A \subseteq A^*(L\beta_{\ell})$;
\item[(iv)] $L\beta_{\ell}$ is a semaphore code.
\ei
\ep

\proof
(i) $\Rw$ (iii). Let $u \in L\beta_{\ell}$ and $a \in A$. Since
$A^*(L\beta_{\ell}) = L \supset A^k$, we may assume that $|u| <
k-1$. Let $b \in A \setminus \{ a\}$ and write $m = k-|u|$. Then
$(a^mu,b^mu) \in \tau_L$, hence
$$(a^{m-1}ua,b^{m-1}ua) = (a^mu \circ a,b^mu \circ a) \in \tau_L.$$
It follows that $a^{m-1}ua$ and $b^{m-1}ua$ must share a suffix in $L$, 
and so $ua$ itself must have a suffix in $L$. Thus
$$(L\beta_{\ell})A \subseteq A^*L = L = A^*(L\beta_{\ell}).$$  

(iii) $\Rw$ (ii). We have
$$LA = A^*(L\beta_{\ell})A \subseteq A^*(L\beta_{\ell}) = L.$$
It follows that $LA^*\subseteq L$. Since $L \in \L_k(A)$, we get $L
\in \I_k(A)$. 

(ii) $\Rw$ (i). By Lemma \ref{irco}, $\tau_L$ is an equivalence
relation.
Let $u,v \in A^k$ be such that $u\tau_L v$. Then $w \leq_s u,v$
for some $w \in 
L$. We may assume that $|w| < k$. Let $a \in A$. Since $L \unlhd
A^*$, we have $wa \in L$. Since $|w| < k$, it follows that $wa$ is a common
suffix of $u\circ a$ and $v\circ a$. Therefore $(u\circ
a)\tau_L(v\circ a)$ and we are done.

(iii) $\iff$ (iv). This follows from Lemma~\ref{alch}, since $L\beta_{\ell}$ is always
a suffix code.
\qed

Note that we can easily produce examples of $L \in \L_k(A) \setminus \I_k(A)$:

\be
\label{lnrc}
Let $A = \{ a,b\}$, $k = 3$ and $L = A^*b \cup A^+Aa$. Then $L \in
\L_k(A)$ but $\tau_L \notin {\rm RC}(A^k)$.
\ee

Indeed, $b \in L$ but $ba \notin L$, hence $L \notin \I_k(A)$ and so
$\tau_L \notin \rc(A^k)$ by Proposition \ref{nsrc}. Note that in this
case $\beta_{\ell} = \{ b,a^3,ba^2, aba, b^2a \}$.

\medskip

Inclusion among left ideals determines inclusion for the equivalence
relations $\tau_L$:

\bl
\label{ordi}
Let $|A| > 1$ and $L,L' \in \L_k(A)$. Then 
$$\tau_L \subseteq \tau_{L'} \iff L \subseteq L'.$$
\el

\proof
Assume that $L \subseteq L'$. Let $(u,v) \in \tau_L$. Then $u$ and $v$
share a common suffix in $L$ and therefore in $L'$. Thus $(u,v) \in \tau_{L'}$.

Assume now that $L \not\subseteq L'$. Let $w \in L \setminus L'$ have
minimum length. Since $A^k \subseteq L'$, we have $|w| < k$. Let $n =
k - |w|$. Fix $a,b \in A$ distinct and take $(u,v) = (a^nw,b^nw) \in A^k \times A^k$. Since
$w \in L$, we have $(u,v) \in \tau_L$. Now $w$ is the longest common
suffix of $u$ and $v$. Since $w \notin L'$, it follows that $(u,v) \notin
\tau_{L'}$. 
\qed

Note that Lemma \ref{ordi} does not hold for $|A| = 1$, since $|A^k| = 1$.

\begin{Definition}
\label{definition.src}
We say that $\rho \in \rc(A^k)$ is a \defn{special right congruence} on $A^k$ if
$\rho = \tau_I$ for some $I \in \I_k(A)$. In view of Proposition
\ref{nsrc}, this is equivalent to say that $\rho = \tau_{A^*S}$ for some
semaphore code $S$ on $A$ such that $A^k \subset A^*S$. We denote by  
$\src(A^k)$ the set of all special right congruences on $A^k$.
\end{Definition}

Note that not every semaphore code $S$ satisfies the condition $A^k \subset A^*S$. However,
it is easy to derive a semaphore code from $S$ that does by considering
\begin{equation}
	S' = (S \cap A^{\le k}) \cup (A^k \setminus A^* S).
\end{equation}
$S'$ is a suffix code since the elements in $S \cap A^{\le k}$ are incomparable in suffix order 
since $S$ is a suffix code, and by construction any element in $A^k \setminus A^* S$ is incomparable 
with the elements in $S \cap A^{\le k}$ and vice versa. Furthermore, $A^k \subset A^*S' \supseteq A^*S$ 
and $SA \subseteq A^* S$ by Lemma~\ref{alch}. Thus $S'A \subseteq A^* S'$ and so by Lemma~\ref{alch} $S'$ 
is a semaphore code.

\bp
\label{isolat}
Let $|A| > 1$. Then:
\bi
\item[(i)] $\tau_{I\cap J} = \tau_I \cap \tau_J$ and $\tau_{I\cup J} =
  \tau_I \cup \tau_J$ for all $I,J \in \I_k(A)$;
\item[(ii)] ${\rm SRC}(A^k)$ is a full sublattice of ${\rm RC}(A^k)$;
\item[(iii)] the mapping
$$\begin{array}{rcl}
\I_k(A)&\to&{\rm SRC}(A^k)\\
I&\mapsto&\tau_I
\end{array}$$
is a lattice isomorphism.
\ei
\ep

\proof
(i) By Lemma \ref{ordi}, we have $\tau_{I\cap J} \subseteq \tau_I \cap
\tau_J$ and $\tau_I \cup \tau_J \subseteq \tau_{I\cup J}$. 

Let $(u,v) \in \tau_I \cap
\tau_J$. Then there exist $x \in I$ and $y \in J$ such that $x \leq_s
u,v$ and $y \leq_s
u,v$. Since $x$ and $y$ are both suffixes of the same word, one of
them is a suffix of the other, say $x \leq_s y$. Then $y \in I \cap J$
and so $(u,v) \in \tau_{I\cap J}$. Thus $\tau_{I\cap J} = \tau_I \cap
\tau_J$.

Assume now that $(u,v) \in \tau_{I\cup J}$. Then there exists some $x
\in I \cup J$ such that $x \leq_s u,v$. If $x \in I$, then $(u,v) \in
\tau_{I}$, otherwise $(u,v) \in \tau_{J}$. Therefore $\tau_{I\cup J} =
  \tau_I \cup \tau_J$.

(ii) Let $I,J \in \I_k(A)$. By part (i), $\tau_{I\cap J}$ is the meet
of $\tau_I$ and $\tau_J$ in both $\rc(A^k)$ and $\src(A^k)$. And 
$\tau_{I\cup J}$ is the join
of $\tau_I$ and $\tau_J$ in both $\rc(A^k)$ and $\src(A^k)$.

Finally, $\tau_{A^kA^*}$ is the identity relation and therefore the
bottom element of both lattices. And 
$\tau_{A^*}$ is the universal relation and therefore the
top element of both lattices.

(iii) This follows from Lemma~\ref{ordi}.
\qed

Given $\rho \in \rc(A^k)$ and $C \in A^k/\rho$, we denote by $\lcs(C)$ the
longest common suffix of all words in $C$. We define
\begin{equation}
\label{equation.rho}
\Lambda_{\rho} = \{ \lcs(C) \mid C \in A^k/\rho \} \qquad \text{and} \qquad
\Lambda'_{\rho} = \{ \lcs(u,v) \mid (u,v) \in \rho \}.
\end{equation}

\bl
\label{propsl}
Let $\rho \in {\rm RC}(A^k)$. Then $A^*\Lambda_{\rho} =
A^*\Lambda'_{\rho} \in \I_k(A)$. 
\el

\proof
Let $C \in A^k/\rho$ and let $w = \lcs(C)$. If $|w| = k$, then $w =
\lcs(w,w)$. If $|w| < k$, then by maximality of $w$ there exist $a,b
\in A$ distinct and $u,v \in A^*$ such that $uaw,vbw \in C$. Thus $w =
\lcs(uaw,vbw)$ and so 
\beq
\label{propsl1}
\Lambda_{\rho} \subseteq
\Lambda'_{\rho}.
\eeq
Therefore $A^*\Lambda_{\rho} \subseteq
A^*\Lambda'_{\rho}$.

Conversely, let $(u,v) \in \rho$. Then $\lcs(u\rho)$ is a suffix of
$\lcs(u,v)$, hence $\Lambda'_{\rho} \subseteq
A^*\Lambda_{\rho}$ and so $A^*\Lambda_{\rho} =
A^*\Lambda'_{\rho}$.

Clearly, $A^*\Lambda'_{\rho} \unlhd_{\ell} A^*$. Since $u = \lcs(u,u)$
for every $u \in A^k$, we have $A^k \subseteq \Lambda'_{\rho}$. Hence
it suffices to show that $(\Lambda'_{\rho})A \subseteq A^*\Lambda'_{\rho}$.

Let $(u,v) \in \rho$ and $a \in A$. We must show that $(\lcs(u,v))a
\in A^*\Lambda'_{\rho}$. Since $A^k \subseteq
\Lambda'_{\rho}$, we may assume that $|\lcs(u,v)| < k-1$. Then
$(\lcs(u,v))a = \lcs(u \circ a,v \circ a)$. Since $(u \circ a,v \circ
a) \in \rho$, we get $(\lcs(u,v))a
\in \Lambda'_{\rho}$ and we are done.
\qed

Given $\rho \in \rc(A^k)$, we write
$$\res(\rho) = \res(\cay(\rho)).$$
We refer to the elements of $\res(\rho)$ as the \defn{resets} of $\rho$.

\bl
\label{props}
Let $\rho \in {\rm RC}(A^k)$. Then:
\bi
\item[(i)] $\res(\rho) = \{ w \in A^* \mid u\rho v \mbox{ for all }u,v \in
A^k \cap(A^*w) \}$; 
\item[(ii)] $\res(\rho) \in \I_k(A)$.
\ei
\el

\proof
(i) Let $w \in \res(\rho)$ and suppose that $u = u'w \in A^k$, $v =
v'w \in A^k$. Since $w \in \res(\rho)$, we have paths
$$p \mapright{u'} p' \mapright{w} r,\quad q \mapright{v'} q'
\mapright{w} r$$ in $\cay(\rho)$. It follows from the definition of
$\cay(\rho)$ that $$u\rho = (u'w)\rho = r = (v'w)\rho = v\rho,$$
hence the direct inclusion holds.

To prove the opposite inclusion, we suppose that $w \in A^* \setminus
\res(\rho)$. Then there exist paths 
$$p' \mapright{w} p, \quad q' \mapright{w} q$$
in $\cay(\rho)$ with $p \neq q$.
If $w$ has a suffix $w'$ of length $k$, then every path
labeled by $w$ ends necessarily in $w'\rho$, hence we must have $|w| <
k$. Since $\cay(\rho)$ is strongly connected by Proposition~\ref{isom}, there exist paths
$$p'' \mapright{x} p', \quad q'' \mapright{y} q'$$
in $\cay(\rho)$ with $|xw| = |yw| = k$. But then
$$(xw)\rho = p \neq q = (yw)\rho$$
and we are done.

(ii) It is immediate that $\res(\rho) \unlhd A^*$. Since every path
in $\cay(\rho)$ labeled by $w \in A^k$ ends necessarily in $w\rho$, we
have $A^k \subseteq \res(\rho)$ and so $\res(\rho) \in \I_k(A)$.
\qed

We can now compare a right congruence with a special right congruence:

\bp
\label{inclu}
Let $|A| > 1$, $\rho \in {\rm RC}(A^k)$ and $I \in \I_k(A)$. Then:
\bi
\item[(i)] $\rho \subseteq \tau_I \iff \Lambda_{\rho} \subseteq I \iff
\Lambda'_{\rho} \subseteq I$;
\item[(ii)] $\tau_I \subseteq \rho \iff I \subseteq \res(\rho)$.
\ei
\ep

\proof
(i) Assume that $\rho \subseteq \tau_I$. Let $(u,v) \in \rho$. Then $u$
and $v$ have a common suffix in $I$, hence $\lcs(u,v)$ has a suffix in
$I$ and so $\Lambda'_{\rho} \subseteq A^*I = I$.

By (\ref{propsl1}), $\Lambda'_{\rho} \subseteq I$ implies
$\Lambda_{\rho} \subseteq I$. 

Finally, assume that $\Lambda_{\rho} \subseteq I$. Let $(u,v) \in
\rho$ and write $w = \lcs(u\rho) \in \Lambda_{\rho} \subseteq
I$. Since $w$ is a suffix of both $u$ and $v$, we get $(u,v) \in
\tau_I$. Thus $\rho \subseteq \tau_I$ as required.

(ii) Assume that $\tau_I \subseteq \rho$. Let $w \in I$ and let $u,v
\in A^k \cap(A^*w)$. Since $u,v$ have a common suffix in $I$, we get
$(u,v) \in \tau_I \subseteq \rho$. Thus $w \in \res(\rho)$ by Lemma
\ref{props}(i) and so
$I \subseteq \res(\rho)$. 

Conversely, assume that $I \subseteq \res(\rho)$. Let $(u,v) \in
\tau_I$. Then we may write $u = u'w$, $v = v'w$ with $w \in I
\subseteq \res(\rho)$. Since $u,v \in A^k \cap(A^*w)$, it follows
from Lemma \ref{props}(i) that $(u,v) \in
\rho$ and so $\tau_I \subseteq \rho$.
\qed

We can now prove several equivalent characterizations of special right
congruences:

\bp
\label{src}
Let $|A| > 1$ and $\rho \in {\rm RC}(A^k)$. Then the following conditions are
equivalent:
\bi
\item[(i)] $\rho \in {\rm SRC}(A^k)$;
\item[(ii)] ${\rm lcs}:A^k/\rho \to A^{\leq k}$ is injective and
  $\Lambda_{\rho}$ is a suffix code;
\item[(iii)] $\rho = \tau_{A^*\Lambda_{\rho}}$;
\item[(iv)] $\rho = \tau_{A^*\Lambda'_{\rho}}$;
\item[(v)] $\rho = \tau_{\res(\rho)}$;
\item[(vi)] $\rho = \tau_L^{\sharp}$ for some $L \in \L_k(A)$;
\item[(vii)] $\Lambda_{\rho} \subseteq \res(\rho)$;
\item[(viii)] $\Lambda'_{\rho} \subseteq \res(\rho)$;
\item[(ix)] whenever
\beq
\label{src3}
p \mapright{aw} q,\quad p' \mapright{bw} q,\quad p'' \mapright{w}
r
\eeq
are paths in ${\rm Cay}(\rho)$ with $a,b \in A$ distinct, then $q = r$.
\ei
\ep

\proof
(i) $\Rw$ (ii).
We start by proving that
\beq
\label{src2}
\lcs(u\tau_I) \in I
\eeq
for all $I \in \I_k(A)$ and $u \in A^k$.

Indeed, for every $w \in u\tau_I$,
there exists some $w' \in I$ such that $w' \leq_s u,w$. Let $z$ be the
shortest suffix among the $w'$. Then $z \in I$ and $z \leq_s w$ 
for every $w \in u\tau_I$, hence $z \leq_s \lcs(u\tau_I)$. Since $I
\unlhd A^*$, it follows that $\lcs(u\tau_I) \in I$ and so
(\ref{src2}) holds.
 
Assume that $\rho = \tau_I$ for some $I \in \I_k(A)$.
We prove that
\beq
\label{src1}
\lcs(u\rho) \leq_s \lcs(v\rho) \Rw (u,v) \in \rho
\eeq
holds for all $u,v \in A^k$. Assume that
$\lcs(u\rho) \leq_s \lcs(v\rho)$. Since
$\lcs(u\rho) \leq_s 
u$ and $\lcs(v\rho) \leq_s v$, it follows that
$\lcs(u\rho)$ is 
a suffix of both $u$ and $v$. 
Now (\ref{src2}) yields $\lcs(u\rho) = \lcs(u\tau_I) \in I$ and
so $u,v$ have a common suffix in $I$.
Therefore $(u,v) \in \tau_I = \rho$ and
(\ref{src1}) holds.

Now (ii) follows from (\ref{src1}).

(ii) $\Rw$ (iii). Write $I = A^*\Lambda_{\rho}$. 
If $(u,v) \in \rho$, then $\lcs(u\rho) \in \Lambda_{\rho} \subseteq I$
is a suffix of both $u$ and $v$, hence $(u,v) \in \tau_I$.

Conversely, let $(u,v) \in \tau_I$. Then there exists some $w \in
\Lambda_{\rho}$
such that $w \leq_s u,v$. Suppose that $\lcs(u\rho) \neq w$. Then
$\lcs(u\rho) <_s w$ or $w <_s \lcs(u\rho)$, contradicting
$\Lambda_{\rho}$ being a suffix code. Hence $\lcs(u\rho) =
w$. Similarly, $\lcs(v\rho) = w$. Since $\lcs:A^k/\rho \to
A^{\leq k}$ is injective, we get $u\rho = v\rho$. Thus $\rho = \tau_I$.

(iii) $\iff$ (iv). This follows from Lemma~\ref{propsl}.

(iii) $\Rw$ (vi). Write $L = A^*\Lambda_{\rho}$. By (iii), we have $\tau_L^{\sharp}
= \rho^{\sharp} = \rho$. Since $L \in \L_k(A)$ by Lemma~\ref{propsl},
(vi) holds.

(vi) $\Rw$ (i). Let $I = LA^* \in \I_k(A)$. Since $L \subseteq I$, it
follows from Lemma \ref{ordi} that $\tau_L \subseteq \tau_I$, hence 
$$\rho = \tau_L^{\sharp} \subseteq \tau_I^{\sharp} = \tau_I$$
by Proposition \ref{nsrc}. 

Now assume that $(u,v) \in \tau_I$. Then there exist factorizations $u
= u'w$ and $v = v'w$ with $w \in I$. Write $w = zw'$ with $z \in
L$. Then $(w'u'z,w'v'z) \in \tau_L$ and so
$$(u,v) = (u'w,v'w) = (u'zw',v'zw') = (w'u'z \circ w',w'v'z \circ w')
\in \tau_L^{\sharp} = \rho.$$  
Thus $\tau_I \subseteq \rho$ as required.

(i) $\Rw$ (v). If $\rho = \tau_I$ for some $I \in \I_k(A)$, then 
$I \subseteq \res(\rho)$ by Proposition \ref{inclu}(ii). Since
$\res(\rho) \in \I_k(A)$ by Lemma \ref{props}(ii), then Proposition
\ref{inclu}(ii) also yields 
$$\tau_{\res(\rho)} \subseteq \rho = \tau_I,$$
hence $\res(\rho) \subseteq I$ by Lemma \ref{ordi}. Therefore $I =
\res(\rho)$.

(v) $\Rw$ (vii) $\iff$ (viii). By Lemma~\ref{props}(ii), $\res(\rho)
\in \I_k(A)$. Now we apply Proposition~\ref{inclu}(i).

(viii) $\Rw$ (i). We have $A^*\Lambda'_{\rho},\res(\rho) \in \I_k(A)$
by Lemmas~\ref{propsl} and~\ref{props}(ii).
It follows from Proposition~\ref{inclu} that
$$\tau_{\res(\rho)} \subseteq \rho \subseteq
\tau_{A^*\Lambda'_{\rho}}.$$
Since $\Lambda'_{\rho} \subseteq \res(\rho)$ yields
$A^*\Lambda'_{\rho} \subseteq \res(\rho)$ and therefore
$\tau_{A^*\Lambda'_{\rho}} \subseteq \tau_{\res(\rho)}$ by Lemma
\ref{ordi}, we get $\rho = \tau_{\res(\rho)} \in \src(A^k)$.

(viii) $\Rw$ (ix). Consider the paths in~\eqref{src3}. Since 
$A^k \subseteq \res(\rho)$ by Lemma~\ref{props}(ii), we may assume
that $|w| < k$. Since $\cay(\rho)$ is strongly connected, there exist
paths
$$s \mapright{x} p, \quad s' \mapright{x'} p'$$
such that $xaw,x'bw \in A^k$. Hence 
$$w = \lcs(xaw,x'bw) \in \Lambda'_{\rho} \subseteq \res(\rho)$$
and so $q = r$.

(ix) $\Rw$ (viii). Let $w \in \Lambda'_{\rho}$. Since 
$A^k \subseteq \res(\rho)$ by Lemma \ref{props}(ii), we may assume
that $|w| < k$. Then $w = \lcs(u,v)$ for some distinct
$\rho$-equivalent $u,v \in A^k$. Hence we may write $u = u'aw$ and $v
= v'bw$ with $a,b \in A$ distinct. Since $u\rho = v\rho$, it follows
that there exist in 
$\cay(\rho)$ paths of the form 
$$s \mapright{u'} p \mapright{aw} u\rho,\quad s' \mapright{v'} p'
\mapright{bw} v\rho.$$ 
Now (ix) implies that $w \in \res(\rho)$.
\qed

\begin{Corol}
\label{corollary.src->semaphore}
If $\rho \in {\rm SRC}(A^k)$ with $|A|>1$, then $\Lambda_\rho$ is a semaphore code.
\end{Corol}

\proof
By Proposition~\ref{src}(ii), $\Lambda_\rho$ is a suffix code. Furthermore, by Lemma~\ref{propsl}
we have $A^*\Lambda_{\rho} \in \I_k(A)$, which in turn implies by Proposition~\ref{nsrc} that
$(A^*\Lambda_\rho)\beta_\ell = \Lambda_\rho$ is a semaphore code.
\qed

We can now prove that not all right congruences are special, even for
$|A| = 2$:

\be
\label{notsp}
Let $A = \{ a,b\}$ and let $\rho$ be the equivalence relation on $A^3$
defined by the following partition:
$$\{ a^3,aba,ba^2\} \cup \{ bab,a^2b \} \cup \{ ab^2 \} \cup \{ b^2a
\} \cup \{ b^3 \}.$$
Then $\rho \in {\rm RC}(A^3) \setminus {\rm SRC}(A^3)$.
\ee

Indeed, it is routine to check that $\rho \in \rc(A^3)$. Since
$\lcs(a^3\rho) = a$ and $\lcs((b^2a)\rho) = b^2a$, then $\Lambda_{\rho}$
is not a suffix code and so $\rho \notin \src(A^3)$ by Proposition~\ref{src}.

\medskip

Let $\rho \in \rc(A^k)$ and let 
\begin{equation}
\label{equation.rho low up}
\begin{split}
	\underline{\rho} &= \vee\{ \tau \in \src(A^k) \mid \tau \subseteq \rho \},\\
	\oo{\rho} &= \wedge\{ \tau \in \src(A^k) \mid \tau \supseteq \rho \}.
\end{split}
\end{equation}
By Proposition \ref{isolat}(ii), we have $\underline{\rho}, \oo{\rho}
\in \src(A^*)$.

\bp
\label{mima}
Let $|A| > 1$ and $\rho \in {\rm RC}(A^k)$. Then:
\bi
\item[(i)] $\underline{\rho} = \tau_{\res(\rho)}$;
\item[(ii)] $\oo{\rho} = \tau_{A^*\Lambda_{\rho}} = \tau_{A^*\Lambda'_{\rho}}$.
\ei
\ep

\proof
(i) By Lemma \ref{props}(ii), we have $\res(\rho) \in \I_k(A)$. Now the
claim follows from Proposition~\ref{inclu}(ii).

(ii) Similarly, we have $A^*\Lambda_{\rho} =
A^*\Lambda'_{\rho} \in \I_k(A)$ by Lemma \ref{propsl}, and the
claim follows from Proposition \ref{inclu}(i).
\qed

The next counterexample shows that the pair
$(\underline{\rho},\oo{\rho})$ does not univocally determine $\rho \in
\rc(A^k)$:

\be
\label{cer}
Let $A = \{ a,b\}$ and let $\rho,\rho'$ be the equivalence relations on $A^3$
defined by the following partitions:
$$\{ a^3,aba,ba^2\} \cup \{ bab,a^2b \} \cup \{ ab^2 \} \cup \{ b^2a
\} \cup \{ b^3 \},$$
$$\{ a^3,b^2a,ba^2\} \cup \{ bab,a^2b \} \cup \{ ab^2 \} \cup \{ aba
\} \cup \{ b^3 \}.$$
Then $\rho,\rho' \in {\rm RC}(A^3)$, $\underline{\rho} =
\underline{\rho'}$ and $\oo{\rho} = \oo{\rho'}$.
\ee

Indeed, we claimed in Example~\ref{notsp} that $\rho$ is a right
congruence, and the verification for $\rho'$ is also straightforward.

It is easy to see that
$$\res(\rho) = A^*A^3 \cup \{ a^2,ab\} = \res(\rho'),$$
hence $\underline{\rho} = \underline{\rho'}$ by Proposition~\ref{mima}(i).

Since
$$\Lambda_{\rho} = \{ a,ab, ab^2, b^2a, b^3 \}$$
and
$$\Lambda_{\rho'} = \{ a,ab, ab^2, aba, b^3 \}$$
we obtain
$$A^*\Lambda_{\rho} = A^+ \setminus \{ b,b^2 \} = A^*\Lambda_{\rho'}$$
and Proposition~\ref{mima}(ii) yields $\oo{\rho} = \oo{\rho'}$.

\medskip

This same example shows also that $\oo{\rho}$ does not necessarily
equal or cover $\underline{\rho}$ in $\src(A^k)$. Indeed, in this case we have
$$\res(\rho) = A^*A^3 \cup \{ a^2,ab\} \subset I \subset A^+ \setminus
\{ b,b^2 \} = A^*\Lambda_{\rho}$$ 
for $I = A^*A^3 \cup \{
a^2,ab,ba\} \in \I_k(A)$. By Lemma~\ref{ordi}, we get
$$\underline{\rho} \subset \tau_I \subset \oo{\rho}.$$

\section{Random walks on semaphore codes}
\label{section.random walks}

As we have seen in Proposition~\ref{mima}, semaphore codes approximate right congruences from above and
below in the lattice structure. In this section, we will define \defn{random walks} (or more specifically \defn{Markov chains})
on semaphore codes. The property that makes this possible is that for a semaphore code $S$ associated to the 
alphabet $A$
\begin{equation}
\label{equation.SA}
	SA \subseteq A^* S,
\end{equation}
see Lemma~\ref{alch}. Namely,~\eqref{equation.SA} implies a right action of $A$ on $S$: for $a\in A$ and
$s\in S$, the action $s.a$ is $t$, if $sa = wt$ with $w\in A^*$ and $t\in S$ under~\eqref{equation.SA}.

To turn the action $S \times A \to S$ into a random walk, we impose a \defn{Bernoulli distribution} on $A^*$,
see~\cite[Section 1.11]{BPR:2010}. More precisely, we associate a probability $0\le \pi(a) \le 1$ to each letter 
$a\in A$ such that $\sum_{a\in A} \pi(a) =1$.  The state space of the random walk is $S$. 
Given $s\in S$, with probability $\pi(a)$ we transition to state $s.a$ in one step. This gives rise to the \defn{transition matrix}
$\mathcal{T}$ with entry in row $s$ and column $s'$
\[
	\mathcal{T}_{s,s'} = \sum_{\substack{a \\ \text{with } s'=s.a}} \pi(a).
\]
Since $\sum_a \pi(a)=1$, it follows that the row sums of $\mathcal{T}$ are equal to one, so that $\mathcal{T}$ is
a row stochastic matrix. Taking $\ell$ steps in the random walk is described by the $\ell$-th power of $\mathcal{T}$, 
that is, the probability of going from $s$ to $s'$ in $\ell$ steps is the $(s,s')$-entry $(\mathcal{T}^\ell)_{s,s'}$ in 
$\mathcal{T}^\ell$. Under the Bernoulli distribution, the probability $\pi(a_1 \cdots a_\ell)$ of a word of length $\ell$
is given by the multiplicative formula $\pi(a_1 \cdots a_\ell) = \prod_{i=1}^\ell \pi(a_i)$.

A suffix code $X$ on $A^*$ is \defn{maximal} if it is not properly contained in any other suffix code on $A^*$, that is,
if $X\subseteq Y \subseteq A^*$ and $Y$ is a suffix code, then $Y=X$. Furthermore, $X$ is called \defn{thin}
if there exists an elements $w\in A^*$ such that $A^* w A^* \cap X = \emptyset$. 
By~\cite[Proposition 3.3.10]{BPR:2010}, for a thin maximal suffix code $X$ we have
$\pi(X) = \sum_{x\in X} \pi(x) = 1$ for all positive Bernoulli distributions $\pi$ on $X$. A Bernoulli distribution on $X$ 
is positive if $\pi(x)>0$ for all $x\in X$.
As shown in~\cite[Proposition 3.5.1]{BPR:2010}, semaphore codes $S$ are thin maximal suffix codes, so that
\begin{equation}
\label{equation.S pi}
	\pi(S) = \sum_{s\in S} \pi(s) = 1.
\end{equation}
Hence any positive Bernoulli distribution on semaphore codes yields a probability distribution.

A \defn{stationary distribution} $I = (I_s)_{s\in S}$ is a vector such that $\sum_{s\in S} I_s = 1$ and
$I\mathcal{T} = I$, that is, it is a left eigenvector of the transition matrix with eigenvalue one. In the finite state case, 
by the Perron--Frobenius Theorem, the stationary distribution exists. It is unique if the random walk is irreducible.
See~\cite{LPW:2009} for more details. In our case, we prove next that a stationary distribution exists and give its
explicit form.

\begin{T}
\label{theorem.stationary}
The stationary distribution of the random walk associated to the semaphore code $S$ is given by
\[
 I = (\pi(s))_{s\in S}\,.
\]
\end{T}

\proof
Taking the $s'$-th component of $I \mathcal{T} = I$ reads
\begin{equation}
\label{equation.sum s}
	\sum_{s\in S} \sum_{\substack{a\in A\\ s'=s.a}} \pi(a) \pi(s) = \pi(s').
\end{equation}
Recall that $s.a=s'$ with $a\in A$ and $s,s'\in S$ means that $sa = ws'$ for some $w\in A^*$. In particular, this can 
only hold if $a$ is the last letter of $s'$ and hence fixed by $s'$.

\noindent {\bf Claim:} The set $S' = \{ w \mid sa = ws', s\in S\}$ for fixed $s'\in S$ with $a\in A$ the last letter of $s'$, 
is a thin maximal suffix code.

Indeed, if the claim is true, we have $\sum_{w \in S'} \pi(w) = 1$ by~\cite[Proposition 3.3.10]{BPR:2010}.
Using that $\pi(a) \pi(s) = \pi(w) \pi(s')$ we can hence rewrite~\eqref{equation.sum s}
\[
	\sum_{s\in S} \sum_{\substack{a\in A\\ s'=s.a}} \pi(a) \pi(s) = \pi(s') \sum_{w\in S'} \pi(w)
	= \pi(s')
\]
as desired. It remains to prove the claim.

First assume that $S'$ is not a suffix code. Then there must be two elements $w,w'\in S'$ that are comparable
in suffix order. But then $ws'$ and $w's'$ are comparable in suffix order, contradicting the fact that $S$ is
a suffix code (since after removing the last letter $a$ the result must be in $S$). Next assume that $S'$ is not maximal.
This means there exists $y \in A^*$ such that $S' \subsetneq S' \cup \{y\}$ is a suffix code. But then
$S\cup \{y \widetilde{s}'\}$ is a suffix code, where $\widetilde{s}'$ is obtained from $s'$ by removing the last letter 
$a$, contradicting the maximality of $S$ (recall that all semaphore codes are maximal 
by~\cite[Proposition 3.5.1]{BPR:2010}). Finally assume that $S'$ is not thin. That means that there exists
$w\in A^*$ such that $A^* w A^* \cap S' \neq \emptyset$. In particular $uwv \in S'$ for some $u,v \in A^*$.
Since by construction $S' \widetilde{s}' \subseteq S$, this would imply $uwv\widetilde{s}' \in S$, contradicting the
fact that $S$ is thin.
\qed

Given $A = \{ a_1,\ldots,a_g \}$ and a right congruence $\rho \in \rc(A^k)$, we are interested in the probability for 
nonempty words of length $\ell \le k$ to be resets on $A^k/\rho$. Since $\res(\rho)=\res(\underline{\rho})$ by
Propositions~\ref{src} and~\ref{mima}, we can restrict ourselves to determine the probabilities for resets of words of
given length for $\underline{\rho} \in \src(A^k)$, or equivalently for semaphore codes $\Lambda_{\underline{\rho}}$
by Corollary~\ref{corollary.src->semaphore}.
 
\bt
\label{theorem.length prob}
Let $\rho \in \rc(A^k)$. Then the probability that a word of length $1\le \ell \le k$ is a reset on $A^k/\rho$ is given by 
\begin{equation}
\label{equation.probability}
	P(\ell) = \sum_{\substack{s \in \Lambda_{\underline{\rho}}\\ \ell(s)\le \ell}} \prod_{a\in s} \pi(a)\;,
\end{equation}
where $a\in s$ in the product runs over every letter in $s$ and
$\ell(s)$ is the length of the word (or suffix) $s$. 
%
\et

\proof
As mentioned above, $\res(\rho)=\res(\underline{\rho})$ by Propositions~\ref{src} and~\ref{mima}
and in addition $\Lambda_{\underline{\rho}}$ is a semaphore code.
Define $\mathrm{Res}(\ell) = \{ w \in A^+ \mid \ell(w)=\ell \text{ and $w$ is a reset on $A^k/\rho$}\}
= \res(\rho) \cap A^\ell$. We claim that 
\[
	\mathrm{Res}(\ell) = \{ w \in A^+ \mid \ell(w)=\ell \text{
          and $w$ has a suffix in $\Lambda_{\underline{\rho}}$}\}. 
\]
Since $\Lambda_{\underline{\rho}}$ is a suffix code, each word has precisely one suffix in $\Lambda_{\underline{\rho}}$.
Hence the claim immediately yields the formula for $P(\ell)$ using that a letter $a\in s$ for $s \in \Lambda_{\underline{\rho}}$
occurs with probability $\pi(a)$.

We prove the claim by induction on $\ell$. By Proposition~\ref{src}(vii) we have that $\Lambda_{\underline{\rho}} 
\subseteq \res(\underline{\rho})=\res(\rho)$. Certainly, for $\ell=1$ the only words that are resets are the words/suffixes of 
length 1 in $\Lambda_{\underline{\rho}}$. Now assume that the claim holds for all words of length less than $\ell$.
Since $\Lambda_{\underline{\rho}} \subseteq \res(\rho)$, we deduce that
\[
	\{ w \in A^+ \mid \ell(w)=\ell \text{ and $w$ has a suffix in $\Lambda_{\underline{\rho}}$}\} \subseteq \mathrm{Res}(\ell)\;.
\]

To prove the reverse inclusion let $v = a_{i_\ell} \ldots a_{i_1} \in \res(\ell)$. If $v\in \Lambda_{\underline{\rho}}$, we are done.
If $a_{i_{\ell-1}} \cdots a_{i_1} \in \mathrm{Res}(\ell-1)$, then by induction $v$ has a suffix in $\Lambda_{\underline{\rho}}$.
Hence assume that $a_{i_{\ell-1}} \ldots a_{i_1} \not \in \mathrm{Res}(\ell-1)$ and $v\not \in \Lambda_{\underline{\rho}}$. 
This requires that $a_{i_\ell} \ldots a_{i_2}$ is a reset, so that again by induction $a_{i_\ell} \ldots a_{i_2}$ has a suffix 
$s$ in $\Lambda_{\underline{\rho}}$. Since $\Lambda_{\underline{\rho}}$ is a semaphore code and 
hence $\Lambda_{\underline{\rho}} A \subseteq A^* \Lambda_{\underline{\rho}}$, we have that if 
$s\in \Lambda_{\underline{\rho}}$, then $s a_{i_1}\in A^*\Lambda_{\underline{\rho}}$. In all cases $v$ has a suffix 
in $\Lambda_{\underline{\rho}}$. This concludes the proof of the claim.
\qed

\be
Take the special right congruence $\rho$ given by congruency classes $\{aaa,baa,aba,bba\}$, $\{aab,bab\}$, 
$\{abb\}$, $\{bbb\}$ with corresponding semaphore code $\Lambda_\rho = \{ a,ab,abb,bbb\}$. The probability to have 
a reset for words of length $\ell$ is
\begin{equation*}
\begin{split}
P(1) & = \pi(a)\\
P(2) & = \pi(a) + \pi(a) \pi(b)\\
P(3) & = \pi(a) + \pi(a) \pi(b) + \pi(a) \pi(b)^2 + \pi(b)^3 = \pi(a) + \pi(a) \pi(b) + \pi(b)^2 = \pi(a) + \pi(b) =1,
\end{split}
\end{equation*}
where for $P(3)$ we have used repeatedly that $\pi(a)+\pi(b)=1$.
\ee

\be
Take the semaphore code
\[
	\{aa, aab, aba, abba, babb, aabb, bbab, abab, bbba, aabb, babbb, abbbb, bbbbb\}\;,
\]
which corresponds to a special right congruence, which is easy to check by Proposition~\ref{src}. Then we have
\begin{equation*}
\begin{split}
P(1) & = 0\\
P(2) & = \pi(a)^2\\
P(3) & = \pi(a)^2 + 2 \pi(a)^2 \pi(b)\\
P(4) &= \pi(a)^2 + 2 \pi(a)^2 \pi(b) + 3 \pi(a)^2 \pi(b)^2 + 3 \pi(a) \pi(b)^3 = \pi(a)^2 + 2 \pi(a)^2 \pi(b) + 3 \pi(a) \pi(b)^2\\
         &= \pi(a)^2 + 2 \pi(a) \pi(b) + \pi(a) \pi(b)^2 = \pi(a) + \pi(a) \pi(b) + \pi(a) \pi(b)^2\\
P(5) &= \pi(a) + \pi(a) \pi(b) + \pi(a) \pi(b)^2 + \pi(a)^2 \pi(b)^3 + 2 \pi(a) \pi(b)^4 + \pi(b)^5\\
        &=  \pi(a) + \pi(a) \pi(b) + \pi(a) \pi(b)^2 + \pi(a) \pi(b)^3 + \pi(b)^4\\
        &= \pi(a) + \pi(a) \pi(b) + \pi(a) \pi(b)^2 + \pi(b)^3 = \pi(a) + \pi(a) \pi(b) + \pi(b)^2\\
        &= \pi(a) + \pi(b) = 1\;,
\end{split}
\end{equation*}
where again we repeatedly used that $\pi(a)+\pi(b)=1$.
\ee

The probability $P(\ell)$ to reach a reset in $\ell$ steps is related to the \defn{hitting time} (see~\cite[Chapter 10]{LPW:2009}).
Namely, given a Markov chain with state space $S$, the hitting time $t_R$ of a subset $R \subseteq S$ is the 
first time one of the nodes in $R$ is visited by the chain. We are interested in the hitting time
$t_{\mathrm{Res}(\rho)}$ for $\rho \in \rc(A^k)$. Set
\[
	p(\ell) = P(\ell) - P(\ell-1) = \sum_{\substack{s \in \Lambda_{\underline{\rho}}\\ \ell(s)= \ell}} \prod_{a\in s} \pi(a)\;.
\]
Then
\[
	t_{\mathrm{Res}(\rho)} = \sum_{\ell=1}^k \ell p(\ell).
\]

\medskip

Note that by Definition~\ref{definition.right congruence}, we also have a right action of $A$ on right congruences
$\rho \in \rc(A^k)$, namely $\rho \times A \to \rho$. Hence, as for semaphore codes, we can define a random walk on
$\rho$ by assigning a probability $\pi(a)$ for each $a\in A$. Recall that by its definition 
in~\eqref{equation.rho low up}, $\underline{\rho}$ is a refinement of $\rho$. Let us relate these various
random walks. A step $s.a=t$ for $s,t \in \Lambda_{\underline{\rho}}$ and $a\in A$ in the random walk on 
the semaphore code $\Lambda_{\underline{\rho}}$ is in one-to-one correspondence to a step $c_s.a=c_t$ in the random 
walk on $\underline{\rho} \in \src(A^*)$, where $c_s,c_t \in \underline{\rho}$ are the unique congruences such that
$\mathrm{lcs}(c_s)=s$, $\mathrm{lcs}(c_t)=t$, respectively. Since $\underline{\rho}$ is a refinement of $\rho$,
a step $c_s.a = c_t$ on $\underline{\rho}$ implies a step $c.a=d$ on $\rho$ whenever $c_s \subseteq c$
and $c_t \subseteq d$. In particular, the transition matrix $\mathcal{T}$ for the random walk on the semaphore code 
$\Lambda_{\underline{\rho}}$ satisfies for a fixed $d\in \rho$
\begin{equation}
\label{equation.lumpable}
	\sum_{\substack{t \in \Lambda_{\underline{\rho}}\\ c_t \subseteq d}} \mathcal{T}_{s,t} = 
	\sum_{\substack{t \in \Lambda_{\underline{\rho}}\\ c_t \subseteq d}} \mathcal{T}_{s',t} 
	\qquad \text{for all $s,s'\in \Lambda_{\underline{\rho}}$ such that $c_{s'} \;\rho\; c_s$.}
\end{equation}
This relation is precisely the condition for a Markov chain to be \defn{lumpable}.
Lumpability was first introduced by Kemeny and Snell~\cite{KS:1976}, see also~\cite[Section 2.3.1]{LPW:2009}.
This means that the transition matrix $\mathcal{T}^\rho$ on $\rho$ indexed by right congruence classes $c,d\in \rho$
can be expressed in terms of $\mathcal{T}$ as follows
\[
	\mathcal{T}^\rho_{c,d} = \sum_{\substack{t \in \Lambda_{\underline{\rho}}\\ c_t \subseteq d}} \mathcal{T}_{s,t} \qquad 
	\text{for any $s\in \Lambda_{\underline{\rho}}$ such that $c_s \subseteq c$.}
\]
The theory of lumpability (or projection) then gives us the stationary distribution $I^\rho$ for $\mathcal{T}^\rho$.

\begin{Proposition}
Let $I^\rho = (I^\rho_c)_{c\in \rho}$ be the stationary distribution for $\mathcal{T}^\rho$. Then
\[
	I^\rho_c = \sum_{\substack{s\in \Lambda_{\underline{\rho}}\\ c_s \subseteq c}} \pi(s).
\]
\end{Proposition}

\proof
By lumpability, we have
\[
	I^\rho_c = \sum_{\substack{s\in \Lambda_{\underline{\rho}}\\ c_s \subseteq c}} I_s,
\]
where $I=(I_s)_{s \in \Lambda_{\underline{\rho}}}$ is the stationary distribution of $\mathcal{T}$.
By Theorem~\ref{theorem.stationary} we have $I_s=\pi(s)$.
\qed

\begin{Remark}
We could have derived an expression for $I^\rho$ also directly from the stationary distribution of the delay
de Bruijn random walk by lumping given as
\[
	I^\rho_c = \sum_{x\in c} \pi(x).
\]
\end{Remark}


\end{document}